\documentclass[final]{siamltex}

\usepackage{graphicx,graphics,epsfig,epsf,listings,url,subfigure}
\usepackage{amsmath,amssymb,amsfonts,mathrsfs}
\usepackage{bbm,bm}
\usepackage{enumitem}
\usepackage{multirow}
\usepackage{color}
\usepackage{algpseudocode,algorithm}

\newcommand{\vect}[1]{\left[ \begin{array}{c} #1 \end{array} \right]}

\usepackage[figuresright]{rotating}

\title{High order implicit-explicit general linear methods \\ with optimized stability regions}

\author{Hong Zhang\thanks
{Computational Science Laboratory, Department of Computer Science, Virginia Polytechnic Institute and State University, Blacksburg, VA 24061 (zhang@vt.edu)}
\and
Adrian Sandu\thanks
{Computational Science Laboratory, Department of Computer Science, Virginia Polytechnic Institute and State University, Blacksburg, VA 24061 (sandu@cs.vt.edu)}
\and
Sebastien Blaise\thanks
{Institute of Mechanics, Materials and Civil Engineering, Universit\'e catholique de Louvain, Louvain-la-Neuve, Belgium (sebastien.blaise@uclouvain.be)}
}

\begin{document}
\maketitle

\begin{abstract}
In the numerical solution of partial differential equations using a method-of-lines approach, 
the availability of high order spatial discretization schemes motivates the development of sophisticated high order 
time integration methods. For multiphysics problems with both stiff and non-stiff terms implicit-explicit (IMEX) 
time stepping methods attempt to combine the lower cost advantage of explicit schemes with the favorable stability properties of implicit schemes.
Existing high order IMEX Runge Kutta or linear multistep methods, however, suffer from accuracy or stability reduction.

This work shows that IMEX general linear methods (GLMs) are competitive alternatives to classic IMEX schemes for large problems arising in practice. 
High order IMEX-GLMs are constructed in the framework developed by the authors \cite{Zhang_2014}. 
The stability regions of the new schemes are optimized numerically. The resulting IMEX-GLMs have similar stability properties as IMEX Runge-Kutta methods, 
but they do not suffer from order reduction, and are superior in terms of accuracy and efficiency. 
Numerical experiments with two and three dimensional test problems illustrate the potential of the new schemes to speed up complex applications. 

\end{abstract}



\section{Introduction}
Many problems in science and engineering are modeled by time-dependent systems of equations involving both stiff and nonstiff terms.
Examples include advection-diffusion-reaction equations,  fluid-structure interactions, and Navier-Stokes equations, and arise in application areas such as mechanical and chemical engineering, astrophysics, meteorology and oceanography, and environmental science.

A method-of-lines approach is frequently employed to separate the spatial and temporal terms in the governing partial differential equations. 
After the spatial terms are discretized by techniques such as finite differences, finite volumes ,and finite elements, 
the resulting system of ordinary differential equations (ODEs) is integrated in time. 
Stiffness may result from different time scales involved (e.g., convective versus acoustic waves), from local processes such as chemical reactions, and from grids with complex geometry \cite{Jebens2011}. 

Explicit numerical integration schemes have maximum allowable time steps bounded by the fastest time scales in the system; 
for example, the time steps are restricted by the CFL stability condition. 
Implicit integration schemes can avoid the step size restrictions but require the solution of large nonlinear systems at each step, 
and are therefore computationally expensive.  
It is therefore of considerable interest to construct numerical integration schemes that avoid the time step restrictions while maintaining a high computational efficiency. 
In the implicit-explicit (IMEX) framework computational efficiency is achieved by performing an implicit integration only for the stiff components of the system.

IMEX methods treat the nonstiff term explicitly and the stiff term implicitly, therefore 
attempting to combine the low cost of explicit methods with the favorable stability properties of implicit methods. 
The development of IMEX linear multistep methods and IMEX Runge-Kutta methods has been reported in 
\cite{Ascher_1995,Frank_1997_IMEXstability,Hundsdorfer_2007_IMEX_LMM,Ascher_1997,Boscarino_2009_IMEXuniform,Pareschi_2000,Verwer_2004}. 

High order methods usually yield more accuracy and better efficiency than low order methods. Many modern PDE solvers are able to employ high order spatial discretizations, e.g., by using high degree polynomials in a discontinuous Galerkin (DG) approach. There is a need to develop high order time stepping formulas to be used in conjunction with high order spatial discretizations. This need motivates the current work. 

Existing high order IMEX methods face challenges when applied to practical problems. 
High order IMEX linear multistep methods suffer from a marked reduction of the stability region with increasing order. 
IMEX Runge-Kutta methods are known to suffer from possible order reduction for stiff problems, 
which reduces the efficiency of high order methods to that of low order methods. 
The order reduction of the former could be avoided by incorporating additional order conditions \cite{Boscarino_2013}. 
Some possible remedies for the latter for Runge-Kutta methods have also been proposed in the literature \cite{Calvo_2001}. 
However, these strategies require special treatment of boundaries which may bring in extra computational cost and complexity; in addition, some of them only work for special cases such as linear boundary conditions.  
To the best of our knowledge, there is no effective way for IMEX Runge-Kutta methods to handle the order reduction in a general way.  
Furthermore, the considerable increase in the number of coupling conditions makes the construction of high order methods difficult. 

This work develops and tests new high order time stepping schemes in the  framework of implicit-explicit general linear methods (IMEX-GLMs) that we have recently developed \cite{Zhang_2014,Zhang_2012}.
The GLM family proposed by Butcher \cite{Butcher_1993a} generalizes both Runge-Kutta and linear multistep methods.
The added complexity gives the flexibility to develop methods with better stability and accuracy properties. 
While Runge-Kutta and linear multistep methods are special cases of GLMs, the framework allows for the construction of many other methods as well. 
In \cite{Zhang_2014,Zhang_2012} we have developed second and third-order IMEX-GLM schemes that showed considerable promise.

This study develops fourth and fifth order IMEX-GLMs with optimized stability properties. Numerical experiments confirm that these methods do not suffer from order reduction, and are considerably more efficient than IMEX-RK methods on a suite of problems ranging from two-dimensional Allen-Cahn and Burgers' equations to three-dimensional compressible Euler equations.

The paper is organized as follows. Section \ref{sec:IMEXGLM} reviews the class of general linear methods.
The construction of high order IMEX-GLMs  with desired stability properties is discussed in Section \ref{sec:HO}.
This section first introduces desirable stability properties building upon existing stability theory for Runge-Kutta methods. 
Numerical results are reported in
Section \ref{sec:results}. Conclusions are drawn in Section \ref{sec:conclusions}.

\section{IMEX general linear methods}\label{sec:IMEXGLM}

IMEX time stepping methods are used to solve systems of ordinary differential equations (ODEs) of the form
\begin{equation}
 \label{eqn:ode_imex}
y' = f(t,y) + g(t,y)\,
\quad t_0 \le t \le t_F\,, \quad y(t_0) = y_0 \in \mathbb{R}^d\,, 
\end{equation}
where $f$ is a nonstiff term, and $g$ is a stiff term.  Many systems of partial differential equations (PDEs) solved in the methods of lines framework lead to partitioned ODE systems \eqref{eqn:ode_imex} after semi-discretization in space. The nonstiff and stiff driving physical processes are captured by $f$ and $g$, respectively.

Partitioned and IMEX general linear methods were developed in \cite{Zhang_2014,Zhang_2012}.  An implicit-explicit general linear method 
applied to \eqref{eqn:ode_imex} advances the solution for one step using:
\begin{subequations}
\label{eqn:imex_glm}
\begin{eqnarray}
Y_i  &=&  h  \sum_{j=1}^{i-1} {a}_{i,j} \, f(Y_j) + h  \sum_{j=1}^{i} \widehat{a}_{i,j} \, g(Y_j)   + \sum_{j=1}^r u_{i,j} \,   y_j^{[n-1]} \,, ~~i=1,\dots,s,~ \\
y_i^{[n]}  &=&   h\sum_{j=1}^s  \left(   {b}_{i,j} \, f(Y_j) +  \widehat{b}_{i,j} \, g(Y_j)  \right)  +  \sum_{j=1}^r v_{i,j} \, y_j^{[n-1]}  \,, ~~ i=1,\dots,r\,.
\end{eqnarray}
\end{subequations}
Such a method is denoted IMEX-GLM$(p,q,s,r)$ ($p$,$q$,$s$ and $r$ stand for order, stage order, number of internal stages, and number of external stages, respectively). The implicit and the explicit components share the same abscissa vector $\mathbf{c}$ and the same coefficients $\mathbf{U}$ and  $\mathbf{V}$.
The IMEX-GLM \eqref{eqn:imex_glm} is represented compactly by the Butcher tableau
\begin{equation}
\label{eqn:imex_glm_tableau}
\renewcommand{\arraystretch}{1.5}
  \begin{array}{c|c|c|c}
\mathbf{c} &   \mathbf{A} & \widehat{\mathbf{A}} & \mathbf{U} \\ \hline 
 & \mathbf{B} & \widehat{\mathbf{B}} & \mathbf{V} \\
\end{array}\,.
\end{equation}
To study the method  \eqref{eqn:imex_glm}  in \cite{Zhang_2014,Zhang_2012}  the additively partitioned original system \eqref{eqn:ode_imex} is written in an equivalent component partitioned form:
\begin{subequations}
\label{eqn:ode_additive}
\begin{eqnarray}
\label{eqn:var-splitting}
y &=& x+z\,, \\
\label{eqn:nonstiff}
x' &=& \tilde{f}(x,z) = f(x+z) \,,\\
\label{eqn:stiff}
z' &=& \tilde{g}(x,z) = g(x+z)\,.
\end{eqnarray}
\end{subequations}

The external vector $y_i^{[n-1]}$ is defined as a $p$th-order approximation of linear combinations of derivatives 
\begin{equation}
\label{eqn:GLM-order_assumption}
 y_i^{[n-1]} = \sum_{k=0}^r {q}_{i,k} h^k x^{(k)} (t_{n-1}) +  \sum_{k=0}^r \widehat{q}_{i,k} h^k z^{(k)}(t_{n-1}) + \mathcal{O}(h^{p+1}), \quad i = 1,\dots,r,
\end{equation}
for some real parameters $q_{i,k}$, $i =1,\dots,r$, $k=0,1,\dots,p$. 
Note that in (\ref{eqn:imex_glm}) $x_i^{[n]}$ and $z_i^{[n]}$ need not to be known individually once they are initialized in the first step. 
Only the combined external vector $y_i^{[n]}=x_i^{[n]}+z_i^{[n]}$ is advanced at each step, similar to how regular GLMs proceed. 

To initialize $y_i^{[0]}$ the starting procedure developed in \cite{Zhang_2014} advances the ODE solution by taking $r-1$ steps with a small step size $\tau$ to obtain the solutions $y_0,y^{\rm start}_1, \dots,y^{\rm start}_{r-1}$. The derivative terms are approximated using only the function evaluations at these $r$ points.  The starting value for the external vector $y_i^{[0]}$ is calculated via the formula
\begin{eqnarray*}
\label{eqn:starting_proc}
y_i^{[0]} &=& y_0 + {q}_{i,1} h f(y_0) + \widehat{q}_{i,1} h g(y_0) \\
&& + \sum_{k=2}^{r} \sum_{j=1}^{r} {q}_{i,k} h^k/\tau^{k-1} {d}_{k,j} f\left(y^{\rm start}_j\right)
+ \sum_{k=2}^{r} \sum_{j=1}^{r} \widehat{q}_{i,k} h^k/\tau^{k-1} {d}_{k,j} g\left(y^{\rm start}_k\right).
\end{eqnarray*}
In vector form it can be written as
\begin{equation}
 y^{[0]} = \mathbf{1}_r \otimes y_0 + \tau \left( \mathbf{Q} \mathbf{D} \otimes \mathbf{I}_{d \times d} \right) \left( (\mathbf{R} \otimes \mathbf{I}_{d \times d} )  F^{\rm start} \right)
    + \tau \left( \widehat{\mathbf{Q}} \mathbf{D} \otimes \mathbf{I}_{d \times d} \right) \left( (\mathbf{R} \otimes \mathbf{I}_{d \times d})  G^{\rm start} \right), 
\end{equation}
where $F^{\rm start}$ and $G^{\rm start}$ consist of function values evaluated at the $r$ starting points, 
e.g. $F^{\rm start} = [f\left(y^{\rm start}_0\right),f\left(y^{\rm start}_1\right), \cdots ,  f\left(y^{\rm start}_{r-1} \right)]^T$.

The $r \times r$ coefficient matrices $\mathbf{Q}$, $\mathbf{D}$, and $\mathbf{R}$ are computed as follows.
\begin{enumerate}
\item $\mathbf{Q}$, $\widehat{\mathbf{Q}}$ are determined by the method coefficients $\mathbf{A}$, $\widehat{\mathbf{A}}$ and the abscissa vector $\mathbf{c}$.  These matrices can be computed column-wise via the order conditions \cite{Butcher_1993a}
\begin{eqnarray}
\label{eqn:initial-weights}
{q}_0 =  \mathbf{1}_s, \quad {q}_i = \frac{\mathbf{c}^i}{i!} - \frac{{\mathbf{A}}\, \mathbf{c}^{i-1}}{(i-1)!}\, ; \quad 
 \widehat{q}_0 =  \mathbf{1}_s, \quad \widehat{q}_i = \frac{\mathbf{c}^i}{i!} - \frac{\widehat{\mathbf{A}}\, \mathbf{c}^{i-1}}{(i-1)!} .
\end{eqnarray}

\item Starting with the following approximation  
\begin{equation}
\vect{\tau x'(t_0) \\ \tau^2 x''(t_0)\\ \vdots \\ \tau^r x^{(r)}(t_0)} =  \tau \mathbf{D} \vect{x'(t_0) \\ x'(t_1)\\ \vdots \\ x'(t_{r-1})} + \mathcal{O}(\tau^{r+1}),
\end{equation}
expanding the right hand side in Taylor series, and comparing the coefficients of each term, allows to identify each entry of $\mathbf{D}$. 

\item $\mathbf{R}$ is a diagonal rescaling matrix which has the form
\begin{equation}
 \mathbf{R} = \textnormal{diag}\left(h/\tau, h^2/\tau^2, \dots, h^r/\tau^r \right). 
\end{equation}
\end{enumerate}
Note that this starting procedure enables to compute the initial approximations with a smaller step size $\tau \le h$. 
The initial approximations can be computed with a regular method of choice; the very small time steps ensure accurate initial solutions,
and also circumvent possible numerical stability issues with the auxiliary scheme. 
The starting procedure used for the experiments in this paper employs the IMEX-RK scheme. 
Considering the possible low accuracy caused by order reduction, 
in the starting procedure we use a step size half as large as the step size for the following integration.  
We point out that using the same step size typically works well based on our experience.

\section{Construction of high order IMEX-GLMs}\label{sec:HO}

We now consider the construction oh high order IMEX-GLMs.
The partitioned GLM theory developed in \cite{Zhang_2014} ensures that, if the stage order is high, the IMEX-GLM method 
has the desired order without the need for coupling conditions. One imposes the order and stage order conditions independently on the implicit and on the explicit component GLMs.

The order conditions for constructing arbitrary GLMs are complicated. In this paper we choose the explicit and implicit components from a subclass of GLMs, named diagonally implicit multistage integration methods (DIMSIMs), for which the order conditions are more manageable.
DIMSIMs are a subclass of GLMs characterized by the following properties \cite{Butcher_1993a}:
\begin{enumerate}
 \item $\mathbf{A}$ is lower triangular with the same element $a_{i,i}=\lambda$ on the diagonal;
 \item $\mathbf{V}$ is a rank-1 matrix with the nonzero eigenvalue equal to one to guarantee preconsistency;
 \item The order $p$, stage order $q$, number of external stages $r$, and number of internal stages $s$ are related by $q\in \{ p-1, p\}$ and $r\in \{ s,s+1\}$.
\end{enumerate}
DIMSIMs can be categorized into four types according to \cite{Butcher_1993a}. 
Type 1 or type 2 methods have $a_{i,j} = 0$  for $j \ge i$ and are suitable for a sequential computing environment, while type 2 and type 3 methods 
have $a_{i,j} = 0$  for $j \ne i$ and are suitable for parallel computation. 
Methods of type 1 and 3 are explicit ($a_{i,i} = 0$), while methods of type 2 and 4 are implicit ($a_{i,i} = \lambda \neq 0$) and potentially useful for stiff systems. 

Following \cite{Zhang_2014}  we are particularly interested in DIMSIMs with $p=q=r=s$, $\mathbf{U}= \mathbf{I}_{s \times s}$, and $\mathbf{V}=\mathbf{1}_s\, v^T$, 
where $v^T\, \mathbf{1}_s=1$ \cite{Jackiewicz_2009_book}. 
The order conditions are satisfied if  the coefficient matrix $\mathbf{B}$ is computed from the relation
\begin{equation} \label{eq:dimsim_coeB}
\mathbf{B}=\mathbf{B}_0-\mathbf{A}\mathbf{B}_1-\mathbf{V}\mathbf{B}_2+\mathbf{V}\mathbf{A},
\end{equation}
where the matrices $\mathbf{B}_0$, $\mathbf{B}_1$,$\mathbf{B}_2$ $\in \mathbb{R}^{s\times s}$ have entries
\[
\left(\mathbf{B}_0\right)_{i,j} = \frac{\int_0^{1+c_i}\phi_j(x)dx}{\phi_j(c_j)},
\quad
\left(\mathbf{B}_1\right)_{i,j} =\frac{\phi_j(1+c_i)}{\phi_j(c_j)},
\quad
\left(\mathbf{B}_2\right)_{i,j} =\frac{\int_0^{c_i}\phi_j(x)dx}{\phi_j(c_j)},
\]
and $\phi_i(x)$ are defined by $\phi_i(x)=\prod_{j=1,j\neq i}^s(x-c_j)$ (cf. \cite[Thm. 5.1]{Butcher_1993a},\cite[Thm. 3.2.1]{Jackiewicz_2009_book}).
Therefore to obtain high order DIMSIMs there is no need to solve complex nonlinear systems as one usually does in the construction of Runge-Kutta methods.  

The important challenge that remains in the construction of IMEX-GLM methods is to achieve the desirable stability properties.   
This section first introduces desirable stability properties building upon existing stability theory for Runge-Kutta methods. 
A numerical optimization process used to maximize the IMEX stability regions is then discussed. 
Two new IMEX-DIMSIM methods of orders four and five are presented at the end.

\subsection{Stability considerations}\label{sec:properties}

\paragraph{A-stability, L-stability, and inherited Runge-Kutta stability}

The classical linear stability theory \cite{Hairer_1993} considers the scalar test problem
whose solution decays to zero
\begin{equation}
\label{eqn:test}
y'=\lambda y\,, \quad t \geq 0\,,  \quad Re(\lambda) \le 0\,.
\end{equation}
A numerical method is stable if, when applied to solve the test problem \eqref{eqn:test} for one step of length $h$ it generates a solution of non-increasing size. A GLM $(\mathbf{A},\mathbf{B},\mathbf{U},\mathbf{V})$ \eqref{eqn:imex_glm_tableau} applied to the test problem gives a solution
\begin{equation}
\label{eqn:GLM_stability_matrixr} 
y^{[n+1]} = \mathbf{M}(z)\, y^{[n]}\,, \quad
 \mathbf{M}(z) = \mathbf{V} + z \, \mathbf{B} \, \left(\mathbf{I}_{s \times s} -z\mathbf{A} \right)^{-1}\, \mathbf{U}\,.
\end{equation}
Here $\mathbf{M}(z)$ is the stability matrix and has a corresponding stability function
\begin{equation}
\label{eqn:GLM_stability_function} 
 p(w,z) = \det(w \mathbf{I}_{r \times r} - \mathbf{M}(z)),
\end{equation}
where $w,z \in \mathbb{C}$ and $z=\lambda h$. 

A-stability requires that the method is unconditionally stable independent of the size of the time step $h$, i.e.,
the spectral radius of the stability matrix $\rho(\mathbf{M}(z)) \le 1$ for any $z$.
L-stability further requires that $\rho(\mathbf{M}(z)) \to 0$ when $z \to \infty$ \cite{Hairer_1993}.
L-stable methods damp components of high frequencies and are particularly useful for stiff problems.
Since IMEX-GLM schemes are designed to treat stiff parts of a given problem implicitly, 
we want the implicit component to be L-stable, or at least A-stable. 
Imposing L-stability directly on the GLM coefficients leads to a difficult analysis, with complexity increasing dramatically as the order increases. 

The inherited Runge-Kutta stability property \cite{Wright_2002,Butcher_2003} provides a practical way to achieve L-stability.
This property requires that the stability function \eqref{eqn:GLM_stability_function}  has the form
\begin{equation}
\label{eqn:IRKS}
 p(w,z) = w^{s-1} \, \bigl(w-R(z)\bigr)\,,
\end{equation}
where $R(z)$ is the stability function of a Runge Kutta method of order $p=s$. 
When \eqref{eqn:IRKS} holds the existing L-stability theory for Runge Kutta methods can be applied to GLMs. 
Note that conditions \eqref{eqn:IRKS} lead to additional nonlinear constraints on method coefficients; these constraints need to be solved accurately in practice.

\paragraph{Stability analysis for IMEX-GLMs}
To study the linear stability of IMEX-GLM schemes we consider the following generalized linear test equation \cite{Zhang_2014}
\begin{equation}
\label{eqn:test-equation}
 y'= \xi y + \widehat{\xi} y\,, \quad t \geq 0,   \quad Re(\xi),Re(\widehat{\xi}) \le 0\,.
\end{equation}
This test problem mimics the structure of \eqref{eqn:ode_imex}.
We consider $\xi y$ to be the nonstiff term and $\widehat{\xi} y$ the stiff term, and denote
 $w=h \xi$ and $\widehat{w}=h \widehat{\xi}$.

Applying (\ref{eqn:imex_glm}) to the test equation \eqref{eqn:test-equation} and assuming $\mathbf{I}_{s \times s}- w \mathbf{A} - \widehat{w} \widehat{\mathbf{A}}$ is nonsingular lead to
\begin{displaymath}
 y^{[n]} = \mathbf{M}(w,\widehat{w})\, y^{[n-1]},
\end{displaymath}
where the stability matrix is defined by \cite{Zhang_2014}
\begin{equation}
 \mathbf{M}(w,\widehat{w}) = \mathbf{V} + \left( w\, \mathbf{B}+ \widehat{w}\, \widehat{\mathbf{B}} \right) \left(\mathbf{I}_{s \times s}-w\, \mathbf{A} - \widehat{w}\,\widehat{\mathbf{A}}\right)^{-1}
 \, \mathbf{U}\,.
\end{equation}
Let $S\subset \mathbb{C}$ and $\widehat{S}\subset \mathbb{C}$ be the stability regions of the explicit GLM component and of the implicit GLM component, respectively.
The {\it combined stability region} is defined by \cite{Zhang_2014}
\begin{equation}
 \label{eqn:stab_region}
 \mathcal{C} =
\left\{\, w \in S, \, \widehat{w} \in \widehat{S}~:~\rho\bigl( \mathbf{M}(w,\widehat{w}) \bigr) < 1 \, \right\} \subset S \times \widehat{S} \subset \mathbb{C} \times \mathbb{C}\,.
\end{equation}
For a practical analysis of stability we define a {\it desired stiff stability region}, e.g.,
\[
\widehat{\mathcal{S}}_\alpha  = \{ \widehat{w} \in \widehat{S} \cap \mathbbm{C}^-~:~ |\mbox{Im}(\widehat{w})| < \tan(\alpha)\, |\mbox{Re}(\widehat{w})| \}\,,
\]
and compute numerically the corresponding {\em constrained non-stiff stability region}:
\begin{equation}
\label{eqn:stability-constrained}
\mathcal{S}_\alpha = \left\{ w \in S ~:~ \rho\bigl( \mathbf{M}(w,\widehat{w}) \bigr) < 1\,,~~\forall\, \widehat{w} \in \widehat{\mathcal{S}}_\alpha \right\}\,.
\end{equation}
The IMEX-GLM method is stable if the constrained non-stiff stability region $\mathcal{S}_\alpha$ is non-trivial (has a non-empty interior) and is sufficiently
large for a prescribed (problem-dependent) value of $\alpha$, e.g., $\alpha = \pi/2$.

\subsection{Finding high order IMEX-DIMSIMs with large stability regions}\label{sec:finding}

The implicit component of the IMEX-GLM is constructed first, and the desired L-stability property is imposed
L-stable GLMs existing in the literature can also be used as implicit components in the combined IMEX scheme.
 
L-stability indicates that $\widehat{w}$ in the non-stiff stability definition (\ref{eqn:stability-constrained}) can be any value on the negative half-plane. 
So the constrained region with $\alpha = \pi/2$ is 
\begin{equation*}
\mathcal{S}_{\pi/2} = \left\{ w \in S ~:~ \rho\bigl( \mathbf{M}(w,r e^{i\theta}) \bigr) < 1\,, ~~\forall\, \theta \in \left[-\frac{\pi}{2},\frac{\pi}{2}\right], ~~\forall\, r \in [0, -\infty) \right\}\,.
\end{equation*}

The corresponding explicit component is constructed next based on the following criteria: it shares the coefficients $\mathbf{c},\mathbf{\widehat{U},}\widehat{\mathbf{V}}$ with the implicit component; it satisfies the desired order conditions;
and results in a large constrained stability region \eqref{eqn:stability-constrained}.

According to the order conditions in \cite{Zhang_2014}, $\mathbf{B}$ depends on $\mathbf{A}$ and $\mathbf{c}$. 
Thus the only free parameters in determining the explicit part are the $s(s-1)/2$ elements of matrix $\mathbf{A}$. 
The problem of finding IMEX-DIMSIMs can be regarded as a numerical optimization problem to find the entries of $\mathbf{A}$
such as to maximize the area of  the constrained stability region $\mathcal{S}_{\pi/2}$. 

We discretize the region $\mathcal{S}_{\pi/2}$ using finite sets of points in polar coordinates
\begin{equation*}
\mathcal{S}_{\pi/2} \approx \left\{ w \in S ~:~ \rho\bigl( \mathbf{M}(w,r e^{i\theta}) \bigr) < 1\,, ~~\forall\, \theta \in \Theta_{\rm f} \subset \left[-\frac{\pi}{2},\frac{\pi}{2}\right], ~~\forall\, r \in R_{\rm f} \subset (-\infty,0] \right\}\,.
\end{equation*}
For example, $R_{\rm f}=[0,-10^{-3},-10^{-2},\dots,-10^3]$ and $\Theta_{\rm f}$ are a set of equally spaced points between $-\pi/2$ and $\pi/2$. 

We next determine the boundary $\partial \mathcal{S}_{\pi/2}$ of the constrained stability region. For this we consider the points of intersection of the boundary with vertical lines on the negative half-plane with abscissae $x_k$.
An intersection point $\widetilde{w}_k=(x_k,y_k)$ should satisfy
\begin{equation}
\label{eqn:boundary}
 \max_{r\in R_{\rm f}, \theta\in \Theta_{\rm f}} \rho\bigl( \mathbf{M}(\widetilde{w}_k,r e^{i\theta}) \bigr) = 1 .
\end{equation}
Note that since the stability region is symmetric, we only need to consider the part above the real axis. 

Starting with an initial point on the vertical line, e.g. $x_k+i\, y_*$ where $y_*$ is large enough to make the point outside the stability region,
we apply the bisection Algorithm \ref{algo:intersection} to find the first point $\widetilde{w}=x_k+i\,y$ along the vertical line such that  
\begin{equation}
\label{eqn:boundary_approx}
 \max_{r\in R_{\rm f}, \theta\in \Theta_{\rm f}} \rho\bigl( \mathbf{M}(\widetilde{w},r e^{i\theta}) \bigr) < 1.
\end{equation}
\begin{algorithm}
 \caption{Bisection algorithm for finding the points of intersection}
  \label{algo:intersection}
 \begin{algorithmic}
  \State Initialize $y_{\rm top}\leftarrow y_0$ $y_{\rm bot} \leftarrow 0$
  \While{ $y_{\rm top}-y_{\rm bot}>tol$ }
    \State $y_{\rm mid}=(y_{\rm top}+ y_{\rm bot})/2$
    \If{ $\widetilde{w} \leftarrow c+i\,y_{\rm mid}$ satisfies the condition \eqref{eqn:boundary_approx} }
      \State $y_{\rm bot} = y_{\rm mid}$
    \Else
      \State $y_{\rm top} = y_{\rm mid}$
    \EndIf
  \EndWhile
  \State \textbf{return} $y_{\rm bot}$
 \end{algorithmic}
\end{algorithm}
A similar idea can be used to find the intersection of the stability region and the real axis, 
which is assumed to be the leftmost point of the stability region.  Then we can determine the boundary with the above-mentioned algorithm. 
Algorithm \ref{algo:area} summarizes the procedure to approximate the area of the stability region.
\begin{algorithm}
 \caption{Algorithm for computing the area of constrained stability regions}
 \label{algo:area}
 \begin{algorithmic}[1]
  \State Find the point $x_b$ of intersection of the stability region and the $x$ axis using a bisection strategy similar to Algorithm \ref{algo:intersection} 
  \State Generate $m$ vertical lines with abscissae $x_k$ linearly spaced between $x_b$ and $0$
  \State Find the points of intersection of these lines and the stability region
  \State Approximate the area of the stability region using the trapezoidal method
 \end{algorithmic}
\end{algorithm}

As we can see, the objective function that approximates the area of the stability region is highly nonlinear and computationally expensive, 
especially for the construction of high order methods. 
The optimization problem is in general difficult to solve numerically. 
First we transform the maximization problem to a minimization problem by minimizing the negative of the objective function.  
Then we use the combination of MATLAB genetic algorithm function, \texttt{ga} and MATLAB local minimizer \texttt{fminsearch}.  
We repeatedly apply the two optimization routines one after another using one's result as the starting point of the other. 
Each optimizer is run multiple times until the results converge; each run is initialized with the previous result. 
We terminate the procedure when the result does not change across multiple runs for both optimizers.

\subsection{New IMEX general linear methods}\label{sec:methods}

The construction of DIMSIMs starts with choosing the abscissa vector $\mathbf{c}$ \cite{Butcher_1996}.
A natural choice is a vector of values equally spaced in the interval $[0,1]$. 
For DIMSIMs of order $p$ and stage order $q=p$, 
the last value $\mathbf{c}_s=1$ allows to use the last stage value as the ODE solution at the next time step. 
This advantage also applies to IMEX-DIMSIM. 
Here we choose the common abscissae for the IMEX pairs equally spaced in $[0,1]$, and including $0$ and $1$.
There is no evidence so far that other choices would lead to better schemes.

\subsubsection{A fourth-order IMEX-DIMSIM pair}
We start with the construction of the implicit part of the IMEX pair. Butcher \cite{Butcher_1996} reports
a failed attempt to construct DIMSIMs with inherited Runge-Kutta stability, $p=q=r=s=4$, and $c=[0,1/3,2/3,1]$. 
Surprisingly we succeeded in solving the nonlinear system comes from the stability constraints by using the Mathematica software. 
For the detailed information on the nonlinear system, we refer to \cite{Butcher_1996}. 
The coefficients of the type 2 (implicit) DIMSIM we found are given in Table \ref{tab:order-4-implicit-coefficients}. 
The choice of the diagonal element of $\widehat{A}$ equal to $0.572816062482135$ ensures that the implicit method $L-$stable,
following the classic theory of Runge-Kutta methods \cite{Hairer_1993}. 
We remark that this new implicit DIMSIM method can be used by itself due to its favorable stability properties.  

\begin{figure}
\centering{
\includegraphics[width=0.86\textwidth]{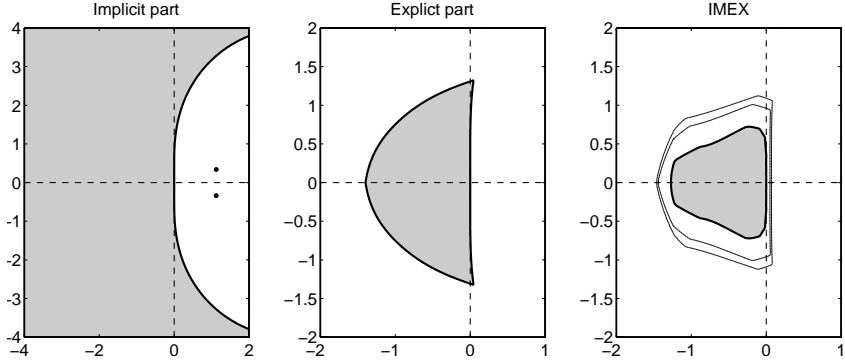}
}
\caption{Stability regions for the fourth-order IMEX-DIMSIM pair with $p=q=r=s=4$ and $c=[0,1/3,2/3,1]$. 
From left to right are stability region $\widehat{S}$ of the implicit method, stability region $S$ of the explicit method, and 
constrained stability regions $\widehat{S}_{\alpha}$ (with $\alpha=\pi/2,\pi/3,\pi/4$ from interior toward exterior, respectively).}
\label{fig:stab_imexdimsim4}
\end{figure}

The optimization problem formulated in Section \ref{sec:finding} for maximizing the constrained stability regions
has six free variables, the lower triangular entries the coefficient matrix $\mathbf{A}$. 
The maximal area of the constrained stability region of the explicit method on the negative plane is approximately $1.34$. 
Figure \ref{fig:stab_imexdimsim4} shows the stability regions of the implicit component $\widehat{S}$, of the explicit component $S$, as well as the constrained stability regions $\widehat{S}_{\alpha}$ for $\alpha=\pi/2,\pi/3,\pi/4$. 

We will refer to the resulting method as IMEX-DIMSIM4. The coefficients of the explicit method to $15$ accurate digits are given in Table \ref{tab:order-4-implicit-coefficients}.

\subsubsection{A fifth-order IMEX-DIMSIM pair}
An L-stable fifth-order type 2 (implicit) DIMSIM with $p=q=r=s=5$ and $c=[0,1/4,1/2,3/4,1]$ was constructed by Butcher \cite{Butcher_1998}. 
We have obtained its coefficients with improved accuracy from $6$ to $15$ decimal digits 
by solving the nonlinear conditions using the Levenberg-Marquardt algorithm implemented by MATLAB's routine \texttt{fsolve}. 

\begin{figure}
\centering{
\includegraphics[width=0.86\textwidth]{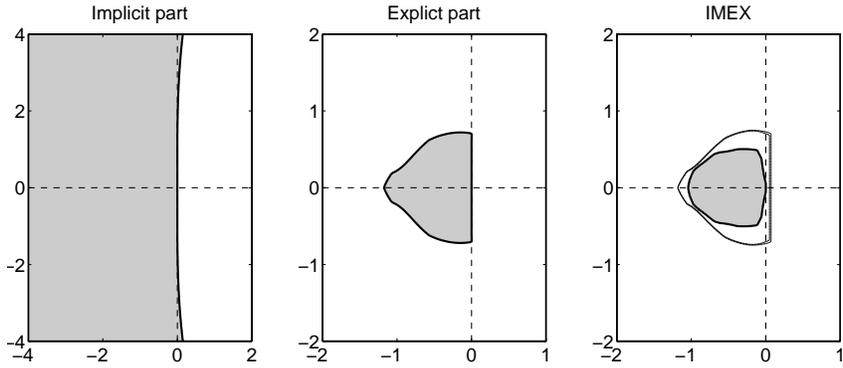}
}
\caption{Stability regions for the fifth-order IMEX-DIMSIM pair with $p=q=r=s=5$ and $c=[0,1/4,1/2,3/4,1]$. From left to right are stability region $\widehat{S}$ of the implicit method, stability region $S$ of the explicit method, and 
constrained stability regions $\widehat{S}_{\alpha}$ (with $\alpha=\pi/2,\pi/3,\pi/4$ from interior toward exterior, respectively) }
\label{fig:stab_imexdimsim5}
\end{figure}

The corresponding explicit component is obtained by the numerical optimization procedure described in the Section \ref{sec:finding}. 
The maximal area of the constrained stability region of the explicit method on the negative plane is approximately $0.83$, and is smaller than the area of the fourth order pair.
Figure \ref{fig:stab_imexdimsim5} shows the stability regions of the implicit component $\widehat{S}$, of the explicit component $S$, as well as the constrained stability regions $\widehat{S}_{\alpha}$ for $\alpha=\pi/2,\pi/3,\pi/4$. 

We will refer to the resulting method as IMEX-DIMSIM5. The coefficients of the method to $15$ accurate digits are given in Table \ref{tab:order-5-implicit-coefficients} (compare the implicit coefficients to \cite{Butcher_1998}). 

\section{Numerical tests}\label{sec:results}
We consider several test problems that are motivated by different application areas such as material science, fluid mechanics, and atmospheric modeling. All problems are governed by partial differential equations and contain both stiff components and nonstiff components.  
The first two test cases are implemented in MATLAB using finite difference schemes for space discretization.  The time integration is performed with the two high order IMEX general linear methods IMEX-DIMSIM4 and IMEX-DIMSIM5. 
The performance of these methods is compared against two classic IMEX Runge-Kutta methods, ARK4(3)6L[2]SA and ARK5(4)8L[2]SA, from Kennedy and Carpenter \cite{Kennedy_2003}. We will refer to these methods as IMEX-RK4 and IMEX-RK5, respectively. 
Both IMEX Runge-Kutta methods have a stiffly-accurate implicit component
and share the same abscissa $\mathbf{c}=\widehat{\mathbf{c}}$ as our IMEX-DIMSIMs do. 

We have also implement the IMEX-DIMSIM schemes in the discontinuous Galerkin solver GMSH-DG \cite{Sebastien_2014} and applied them to the three-dimensional compressible Euler equations coming from multiscale nonhydrostatic atmospheric simulations. 

All the experiments have been performed on a workstation with four Intel Xeon E5-2630 Processors. 
The goal is to assess the performance of the high order IMEX-DIMSIM and IMEX-RK methods on both two-dimensional and three-dimensional simulations.

\subsection{Allen-Cahn equation}\label{sec:allencahn}
We consider the two-dimensional reaction-diffusion Allen-Cahn  problem \cite{Chen_2004} which describes the process of phase transition in materials science.
\begin{equation}
\label{eqn:allen-cahn}
\frac{\partial u}{\partial t} = \alpha \nabla^2 u + \beta (u-u^3) + f, \quad
0 \le x,y \le 1\,, \quad 0 \le t \le 0.5,
\end{equation}
where the parameters are $\alpha=0.01$, $\beta =3.$, and $f(t,x,y)$ is a source term that is consistent with the exact solution
 $u(t,x,y) = 2 + \sin(2 \pi (x-t)) \cos(3\pi(y-t))$. Time varying Dirichlet boundary conditions (that represent the exact solution evaluated at the boundaries) are imposed.   The spatial discretization is performed using a second-order central finite difference scheme on a uniform grid with $\Delta x =\Delta y = 1/40$.  

Explicit time stepping methods have a maximal allowable time step  $h \propto \Delta x ^2 $ due to the CFL condition related to diffusion. 
To overcome this limitation we treat the stiff diffusion term implicitly and the remaining terms explicitly. 
Since the discrete diffusion term is linear we perform a single LU factorization of the matrix $\mathbf{I} - h\gamma \mathbf{J}$ and reuse it throughout the simulation; here $\gamma$ is a method coefficient and $\mathbf{J}$ is the Jacobian of the stiff diffusion.  
  
\begin{figure}
  \subfigure[Convergence diagram]{
  \label{fig:allencahn-conv}
    \centering{
      \includegraphics[width=0.45\textwidth]{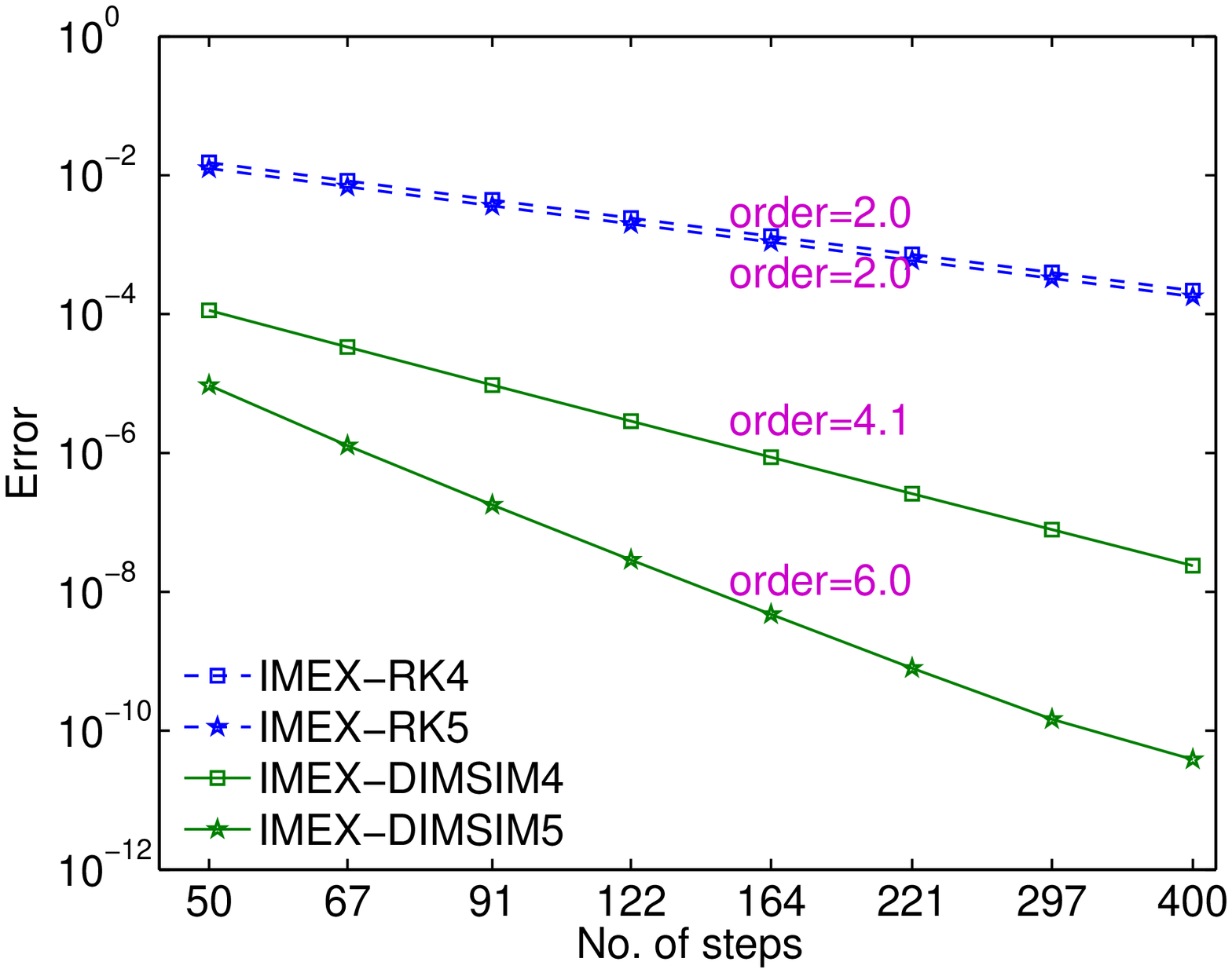}
    }
  }
  \subfigure[Work-precision diagram]{
  \label{fig:allencahn-wp}
  \includegraphics[width=0.45\textwidth]{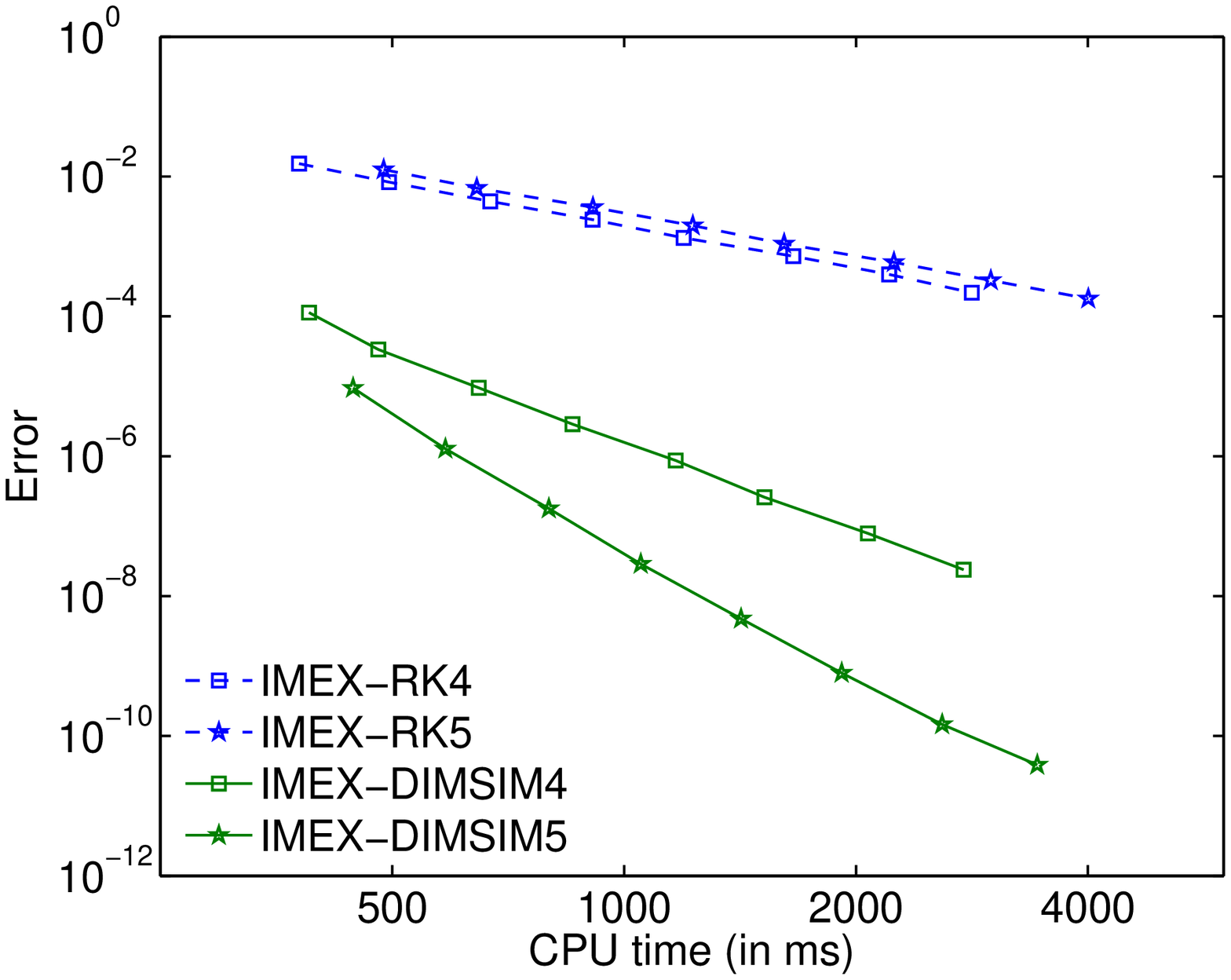}
  }
\caption{Comparison of high order IMEX-DIMSIM and IMEX-RK results for the 2D Allen-Cahn equation \eqref {eqn:allen-cahn}. Shown are the temporal discretization errors corresponding to the solution at the final time $t=0.5$. 
}
\label{fig:allencahn}
\end{figure}

\begin{figure}
    \centering{
      \includegraphics[width=0.8\textwidth]{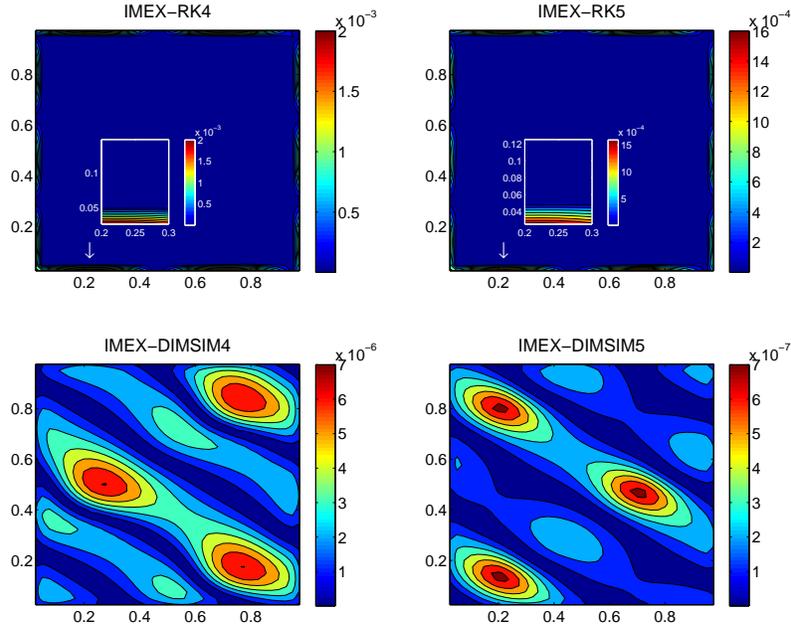}
    }
\caption{Absolute temporal errors at the final time $t=0.5$ for various IMEX schemes on the 2D Allen-Cahn equation  \eqref {eqn:allen-cahn}. 
A fixed time step of size $h=1/50$ is used. IMEX-RK methods show large errors originating near boundaries. IMEX-DIMSIM methods have much smaller errors which are distributed over the entire domain.}
\label{fig:allencahn_error}
\end{figure}

The reference solution $u_{\rm ref}$ is obtained using MATLAB's routine \texttt{ode15s} with very tight tolerances $AbsTol=RelTol=3\times10^{-14}$. The absolute solution error magnitude is measured in the $L_2$ norm:
\begin{equation}
\label{eqn:L2error}
\mathbf{E} = \|u-u^{\rm ref}\|_2. 
\end{equation}
Figure \ref{fig:allencahn-conv}  shows the errors at the final time for solutions computed using different numbers of steps. 
The two IMEX-RK methods show a marked order reduction - to order two. There is no order reduction for the IMEX-DIMSIM schemes; IMEX-DIMSIM4 displays the theoretical order while IMEX-DIMSIM5 shows a higher convergence than the theoretical order.
The IMEX-DIMSIMs give considerably more accurate results than the IMEX-RK methods for all step sizes tested. 
This is noteworthy
since IMEX-DIMSIMs have fewer stages than the IMEX-RK methods of the same order and therefore require fewer function evaluations and linear solves per step.
The corresponding work-precision diagrams of errors versus CPU time are shown in Figure \ref{fig:allencahn-wp} and reveal a sizable
gap in efficiency between the two families of IMEX schemes.
Figure \ref{fig:allencahn_error} shows the spatial distribution of the absolute errors $|u_{\rm numerical}-u_{\rm reference}|$ at final time; this is only the temporal discretization error as we compare against a reference solution that uses the same spatial discretization.  IMEX-RK methods give large errors near boundaries and relatively smaller errors in the interior of the domain are evenly distributed.   The order reduction phenomenon of IMEX-RK methods originates with errors at the boundaries, but plague the whole domain as the time evolves. 
In contrast, IMEX-DIMSIMs handle the boundaries well and preserve the theoretical orders of convergence. 
%

\subsection{Burgers' equation}\label{sec:burger}
The two-dimensional viscous Burgers equation \cite{Bahadir_2003}
\begin{equation}
\label{eqn:burgers}
\frac{\partial u}{\partial t} + \frac{1}{2} \nabla (u \cdot u) = \nu \nabla^2 u, \quad \nu = 0.1, \quad 0 \le x,y \le 1, \quad 0 \leq t \leq  1,
\end{equation}
is a simplification of the 2D Navier-Stokes equations which admits the analytic solution 
\[
u^{\rm analytic}(t,x,y) = \left( 1+e^{\frac{x+y-t}{2\, \nu}}\right)^{-1}.
\]
The initial conditions and the Dirichlet boundary values correspond to the analytic solution. 
Spatial derivatives are discretized with second order central finite differences on a uniform grid 
with resolution $\Delta x =\Delta y = 1/50$. 

The application of the IMEX integration treats the diffusion term implicitly and the convective term explicitly.
We compare the numerical solutions against a reference solution computed with MATLAB routine $ode15s$ with tolerances $AbsTol=RelTol=3\times10^{-14}$ that uses the same spatial discretization.
Therefore the errors \eqref{eqn:L2error} reported here are only due to the temporal discretization.
 
\begin{figure}
  \subfigure[Convergence diagram]{
\label{fig:burgers-conv}
    \centering{
      \includegraphics[width=0.45\textwidth]{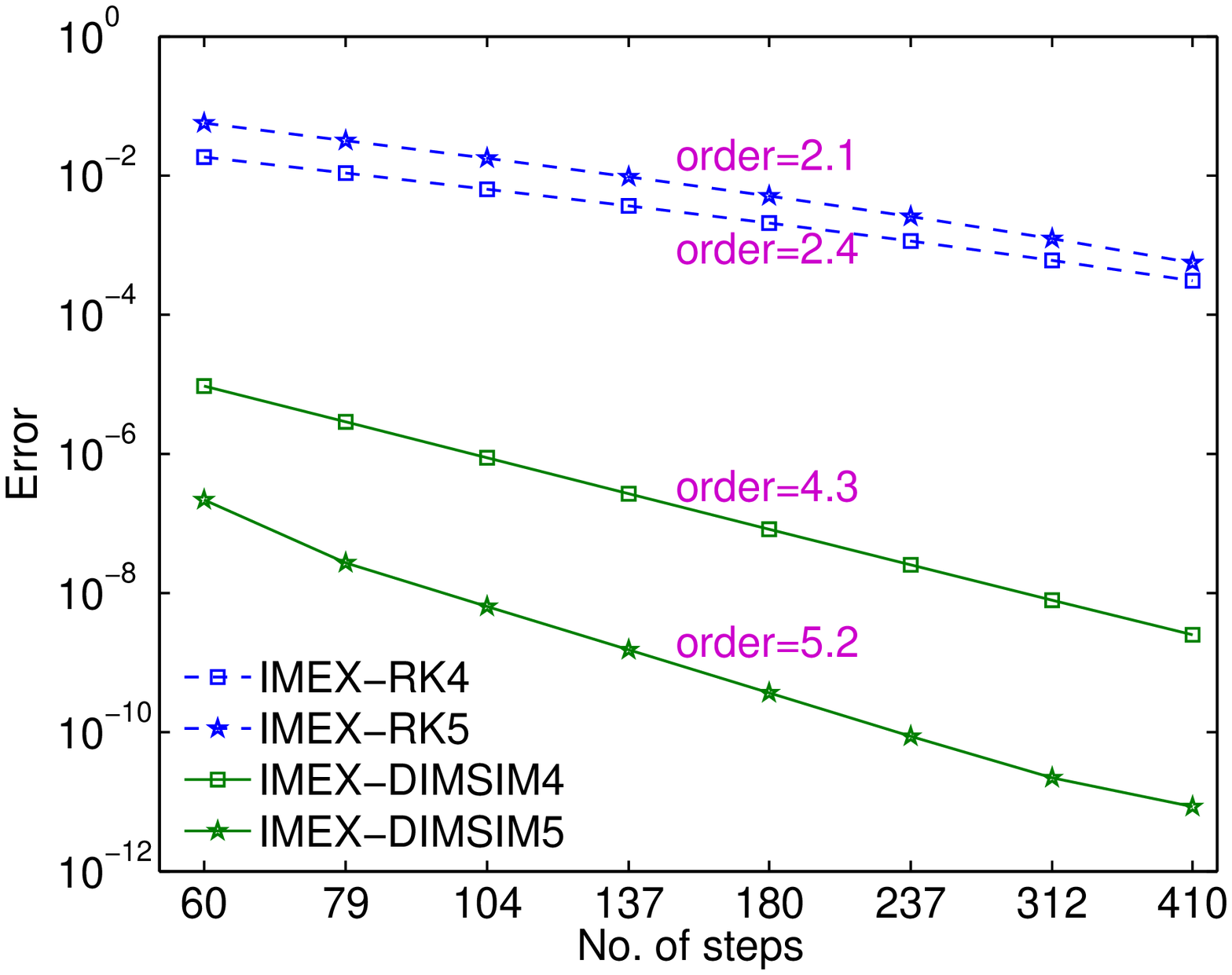}
    }
  }
  \subfigure[Work-precision diagram]{
\label{fig:burgers-wp}
  \centering{
  \includegraphics[width=0.45\textwidth]{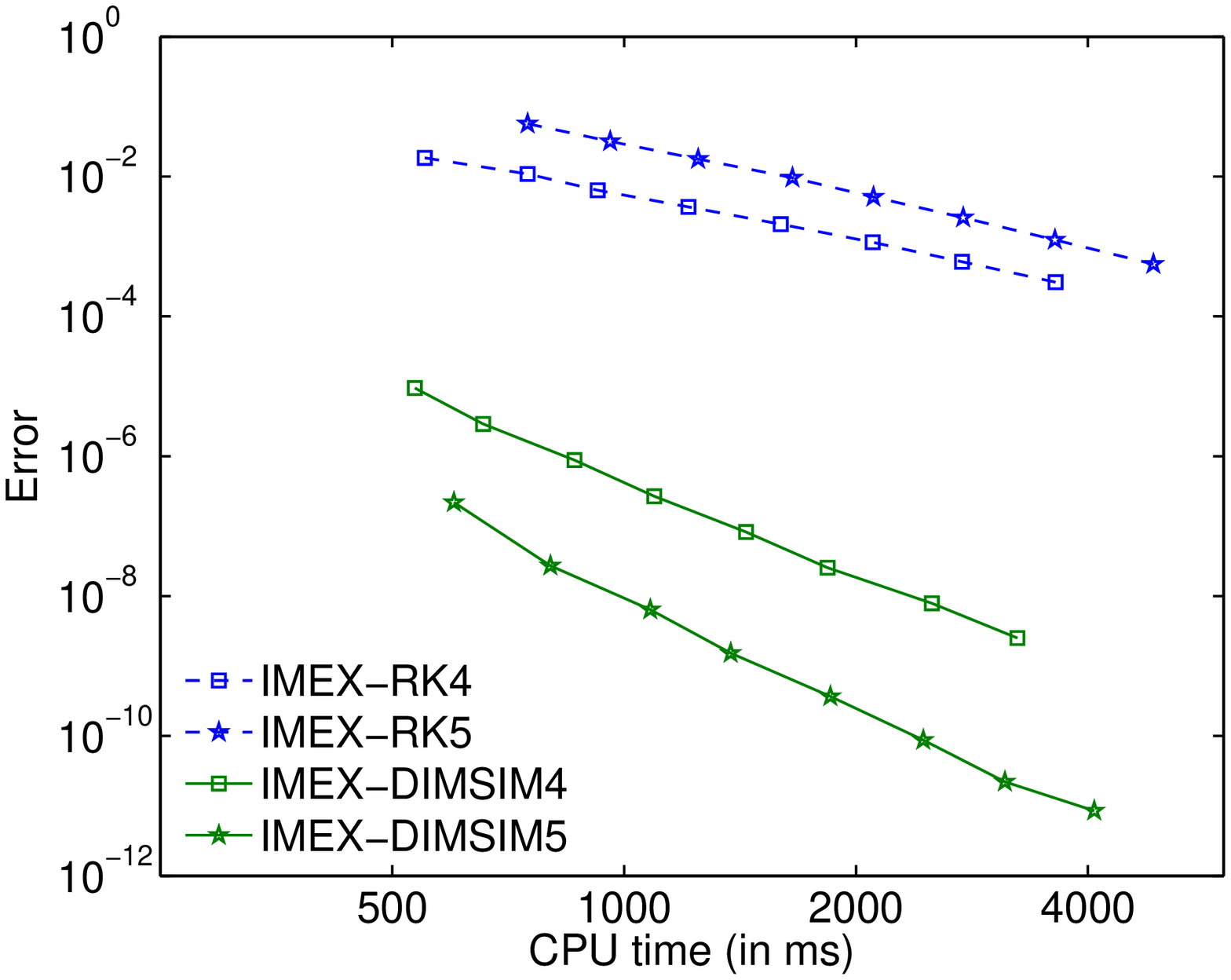}
  }
  }
\caption{Comparison of high order IMEX-DIMSIM and IMEX-RK results for the 2D viscous Burgers equation \eqref{eqn:burgers}. The integration time interval is $[0,1]$. Shown are the temporal discretization errors corresponding to the solution at the final time $t=1$. 
}
\label{fig:burgers}
\end{figure}

Figure \ref{fig:burgers} compares the performance of the high order IMEX schemes. The convergence diagram in Figure \ref{fig:burgers-conv}
reveals that two IMEX-RK methods show order reduction to order two. The two IMEX-DIMSIMs converge with their theoretical orders.
The efficiency diagram in Figure \ref{fig:burgers-wp} illustrates again a gap in performance between the two families, with
IMEX-DIMSIMs demonstrating a considerably better efficiency than IMEX-RK methods. 

\begin{figure}
    \centering{
      \includegraphics[width=0.8\textwidth]{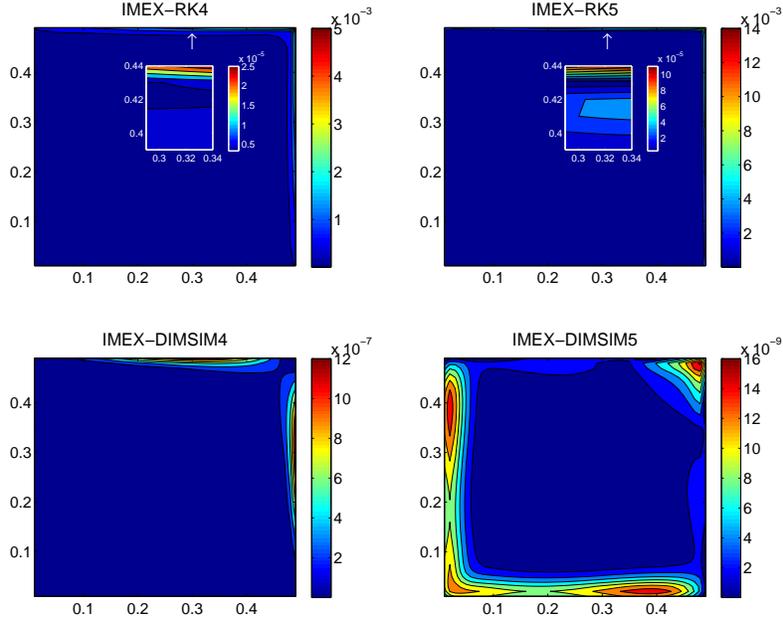}
    }
\caption{Absolute temporal errors at the final time $t=1$ for various IMEX schemes on the 2D viscous Burgers equation \eqref{eqn:burgers}. 
A fixed time step of size $h=1/50$ is used. All methods show larger errors originating near boundaries. IMEX-DIMSIM methods have much smaller errors overall.}
\label{fig:burgers_error}
\end{figure}

Figure \ref{fig:burgers_error} shows the spatial distribution of absolute errors at the final time.
The boundary errors dominate the accuracy of the results for  all schemes. 
The boundary conditions for this PDE may be more challenging than the previous one 
since they affects both spatial derivative terms in \eqref{eqn:burgers}. Nevertheless, the error
magnitude is much smaller for the IMEX-DIMSIM solutions.
%
 
\subsection{Application to atmospheric simulations}\label{sec:rising bubble}
\subsubsection{Compressible Euler equations}
The dynamics of non-hydrostatic atmospheric processes can be described by the compressible Euler equations \cite{Giraldo_2008}:
\begin{subequations}
 \label{eqn:euler}
\begin{eqnarray}
 \nonumber
\frac{\partial \rho}{\partial t} + \nabla \cdot \left( \rho \mathbf{u} \right) & = & 0 \\
 \label{eqn:gw-equation}
 \frac{\partial \rho \mathbf{u}}{\partial t} + \nabla \cdot \left( \rho \mathbf{u} \mathbf{u} + p \mathbf{I}\right) & = & - \rho g \mathbf{\widehat{e}_z} \\
 \nonumber
\frac{\partial \rho \theta}{\partial t} + \nabla \cdot \left( \rho \theta\mathbf{u} \right) & = & 0\,,
 \end{eqnarray}
where  $\rho$ is the density, 
$\mathbf{u}=(�u,v,w)^T$ is the velocity vector, $w$ being used in three-dimensional case, 
$\theta$ is the potential temperature, and $\mathbf{I}$ is the identity matrix. 
The gravitational acceleration is denoted by $g$ while $\mathbf{\widehat{e}_z}$ is a unit vector pointing upwards.
The prognostic variables are $\rho$, $\rho \mathbf{u}$, and $\rho \theta$.
The pressure $p$ in the momentum equation is computed by the equation of state
\begin{equation}
 \label{eqn:euler-pressure}
 p=p_0\left( \frac{\rho \theta R_d}{p_0} \right)^{\frac{c_p}{c_v}},
\end{equation}
\end{subequations}
where $p_0=10^5$ Pa is the surface pressure, $R_d$ is the ideal gas constant, and $c_p$ and $c_v$ are the specific heat of the air for constant pressure and volume. 
To better maintain the hydrostatic state we follow the splitting introduced by Giraldo and Restelli \cite{Giraldo_2008}
\begin{eqnarray*}
\rho(\mathbf{x},t) &=& \bar{\rho}(z) + \rho'(\mathbf{x},t) \\
(\rho \theta)(\mathbf{x},t)& =& \overline{(\rho \theta)} (z) + (\rho \theta)'(\mathbf{x},t) \\
p(\mathbf{x},t) &=& \bar{p}(z) + p'(\mathbf{x},t),
\end{eqnarray*}
where the overlined values are in hydrostatic balance. 
The governing equation \eqref{eqn:euler} can then be rewritten as
\begin{subequations}
\label{eqn:new-euler}
\begin{eqnarray}
\nonumber
 \frac{\partial \rho'}{\partial t} & = & - \nabla \cdot \left( \rho \mathbf{u} \right)  \\
\label{eqn:new-gravity-wave-eqn}
 \frac{\partial \rho \mathbf{u}}{\partial t} & = & - \nabla \cdot \left( \rho \mathbf{u} \mathbf{u} + p' \mathbf{I}\right) - \rho' g \mathbf{\widehat{e}_z} \\
\nonumber
 \frac{\partial (\rho \theta)'}{\partial t} & = & - \nabla \cdot \left( \rho \theta\mathbf{u} \right) \,,
\end{eqnarray}
and closed with
\begin{equation}
 p'=p_0\left( \frac{\rho \theta R_d}{p_0} \right)^{\frac{c_p}{c_v}}- \bar{p}. 
\end{equation}
\end{subequations}
The equations are discretized in space using the discontinuous Galerkin method, whose usage for geophysical simulations is gaining popularity, e.g. \cite{blaise2010,Blaise_2012,nair2005,stcyr2009,Giraldo_2008}. 
The model, based upon the mesh database of the GMSH mesh generator code \cite{Gmsh}, has been used to solve several PDEs, either in the domain of geophysics \cite{Seny_2013,karna2013} and engineering \cite{seny2014,kameni2012}.
For more information about the space discretization, refer to \cite{Sebastien_2014}.

The set of equations \eqref{eqn:new-euler} applied to atmospheric flows is a good candidate for an IMEX time discretization, because of the different temporal scales involved. 
In usual atmospheric configurations, the acoustic waves are the fastest phenomena, with a propagation speed of about $340$ ms$^{-1}$. 
This high celerity restricts the explicit time step to a small value due to the CFL stability condition. 
However, acoustic waves are generally not important for the modeler who is more interested by advective timescales. 
The IMEX method allows to circumvent the CFL condition by treating the linear acoustic waves implicitly, while the remaining terms are explicit. 
According to Giraldo et al. \cite{Giraldo_2010}, the right-hand side of \eqref{eqn:new-gravity-wave-eqn} is additively split into a linear part responsible for the acoustic waves and a nonlinear part. 
The linear term
\begin{equation}
 -
 \begin{bmatrix}
  \nabla \cdot \left( \rho \mathbf{u} \right) \\
  \nabla \cdot \left( p' \mathbf{I}\right)  + \rho' g \mathbf{\widehat{e}_z}\\
  \nabla \cdot \left( \rho \bar{\theta}\mathbf{u} \right)
 \end{bmatrix}
\end{equation}
with the pressure linearized as 
\begin{equation*}
 p'=\frac{c_p \bar{p}}{c_v\overline{\rho \theta}} \left( \rho \theta \right)'
\end{equation*}
is treated implicitly, while the remaining (nonlinear) terms are treated explicitly. 

\subsubsection{Test cases}

In this paper we consider two-dimensional and three-dimensional rising thermal bubble test cases slightly modified from the ones introduced in \cite{Giraldo_2008}. 
\paragraph{Two-dimensional case} 
The motion of the air is driven by a time varying potential temperature perturbation from the bottom boundary
\begin{equation}\label{eqn:bottom}
\theta' = \left\{
\begin{array}{l l}
 0 & \textnormal{for} ~~ r >r_c, \\
 \frac{\theta_c}{2} \left(1+\cos \left(\frac{\pi r}{r_c}\right) \right) \sin^2 \left(\frac{\pi t}{50} \right) & \textnormal{for} ~~ r \leq r_c,
\end{array}
\right.
\end{equation}
where $\theta_c=5 ^\circ C$, $r=\sqrt{(x-x_c)^2 }$, $r_c=250$ m, and $(x,z) \in [0,1000]^2$ with $t\in [0,200]$ s
and $x_c=500$ m. 
No-flux boundaries are used for the other three boundaries. 
The computational domain is a 2D uniform mesh with actual resolution of about $66.7\times66.7$ m.
Fourth-order polynomials are used on each element. 
The resulting ODE system contains 
$\sim 2.3 \times 10^{4}$ variables.

\paragraph{Three-dimensional case} 
Diffusion terms 
\begin{equation}
 \begin{bmatrix}
  \nabla \cdot \left( \mu \nabla \rho'\right) \\
  \nabla \cdot \left( \mu \nabla (\rho \mathbf{u}\right)) \\
  \nabla \cdot \left( \mu \nabla (\rho\theta)' \right)
 \end{bmatrix}
\end{equation}
with $\mu=6$ m$^2$s$^{-1}$ are added to the right-hand side of \eqref{eqn:new-gravity-wave-eqn} to limit the oscillations resulting from a high order spatial discretization of a complex flow on a coarse grid.

The bottom boundary is also imposed as \eqref{eqn:bottom} with $r=\sqrt{(x-x_c)^2 + (y-y_c)^2}$,  $r_c=250$m, $(x,y,z) \in [200,800]^2 \times[0,600]$  and $(x_c,y_c)=(500,500)$. 
No-flux boundaries are used for all the other boundaries. 
Considering the more expensive computational cost of the 3D test, we use a polynomial order of $3$ for the DG scheme. 
A 3D uniform mesh grid with actual resolution of  $100\times100\times100$ m is used. 
The resulting ODE system has 
$\sim 7 \times 10^{4}$ degrees of freedom.

\begin{figure}
  \subfigure[2D solution at final time $t=200s$.]{
  \label{fig:snapshots2}
    \centering{
      \includegraphics[width=0.45\textwidth]{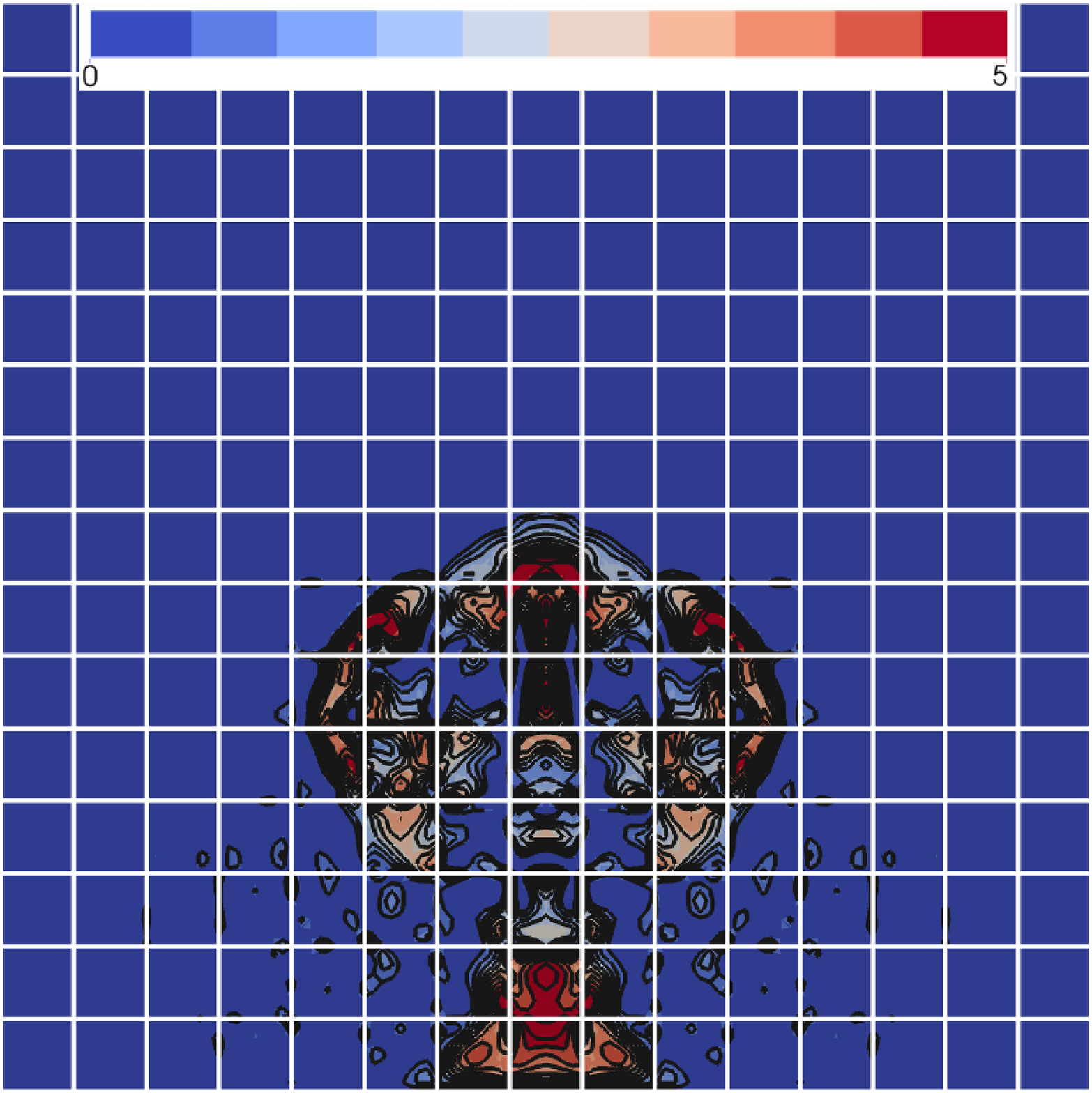}
    }
  }
  \subfigure[3D solution at final time $t=300s$.]{
  \label{fig:snapshots3}
  \centering{
  \includegraphics[width=0.45\textwidth]{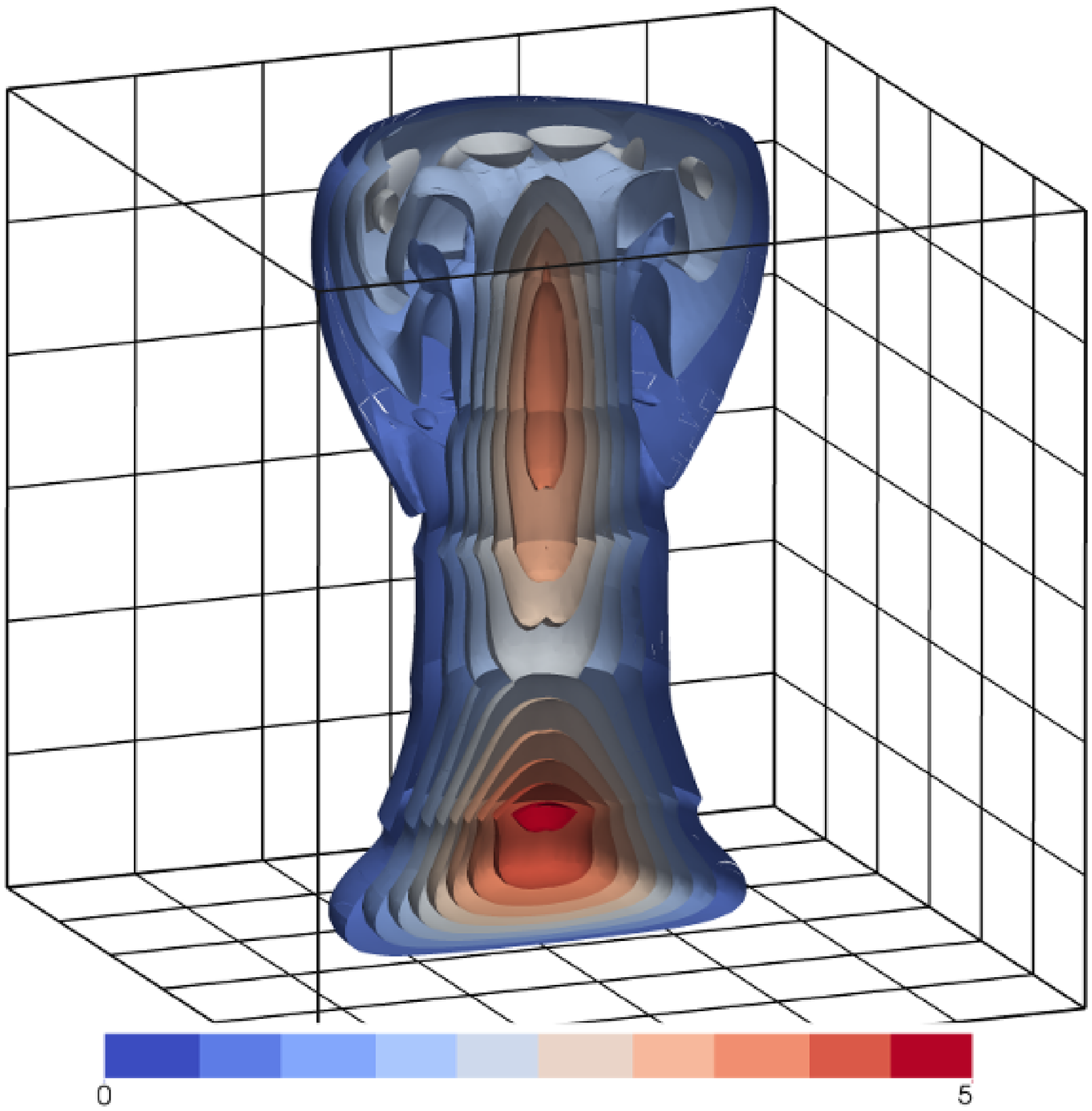}
  }
  }
\caption{Perturbation of potential temperature (in $^\circ C$) from the simulation of thermal rising bubble \eqref{eqn:new-euler}. The background mesh is displayed in wireframe.}
\label{fig:snapshots}
\end{figure}

Figure \ref{fig:snapshots} shows the reference solutions at the final time for 2D and 3D cases. 
%

\subsubsection{Numerical results}
The relative $L_2$ errors for each of the prognostic variables 
\begin{equation}
\label{eqn:L2err}
 \mathbf{E}(q)=\sqrt{\frac{\int_\Omega (q^{\rm numerical}-q^{\rm reference})^2 d\Omega}{ \int_\Omega (q^{\rm reference})^2 d\Omega}},
\end{equation}
are measured against a reference solution obtained by applying the classic fourth-order explicit RK method to solve the original (non-split) model with a very small time step $h=0.005s$. 
Since the time varying boundary conditions are imposed directly on the temperature term $p'$ in the momentum equations of \eqref{eqn:new-euler}, we discuss the results for the variables $\rho \mathbf{u}$.

\begin{figure}    \centering{
  \subfigure[Convergence diagram]{

      \includegraphics[width=0.66\textwidth]{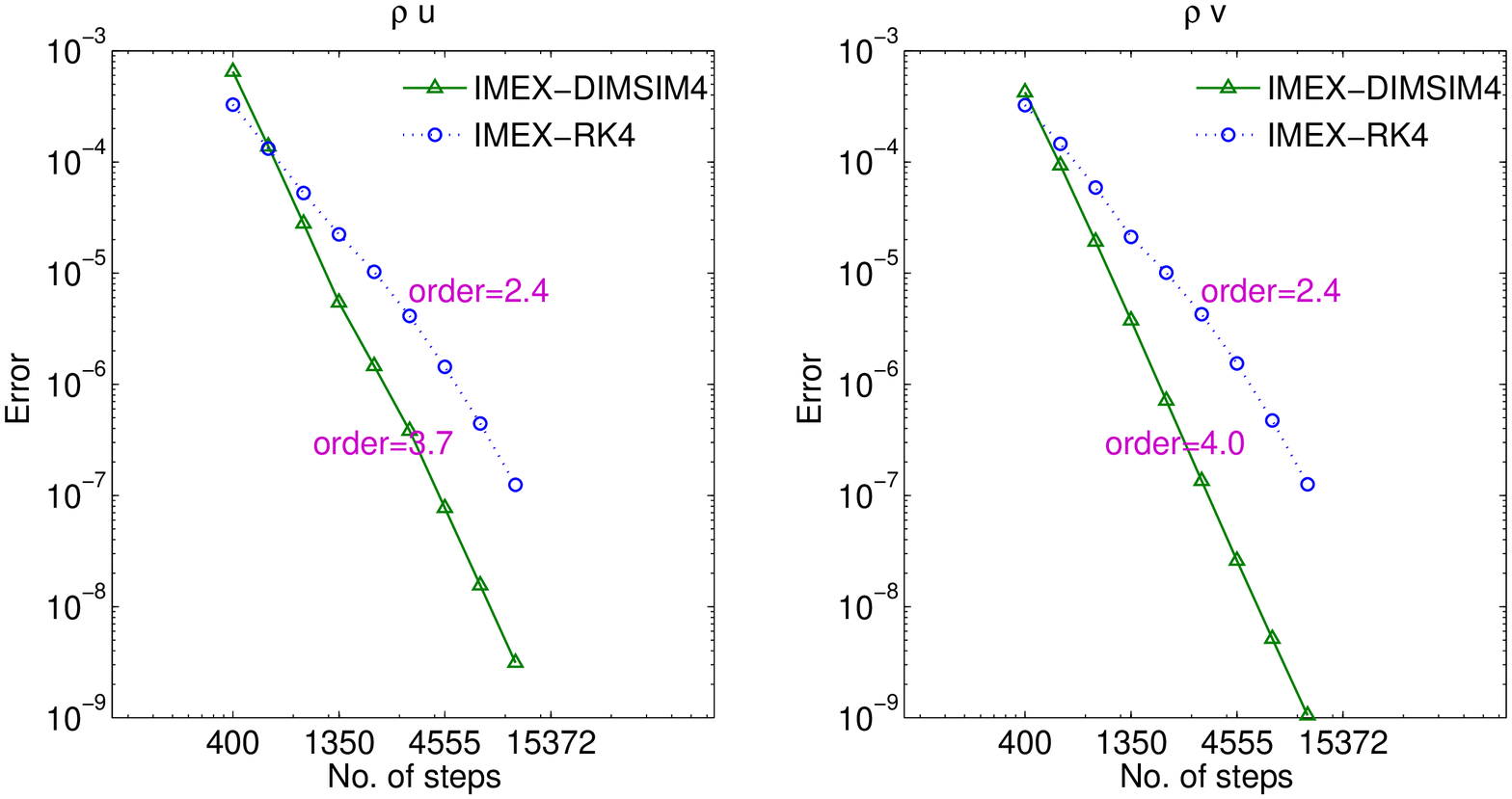}
 
  }
  \subfigure[Work-precision diagram]{
  \includegraphics[width=0.66\textwidth]{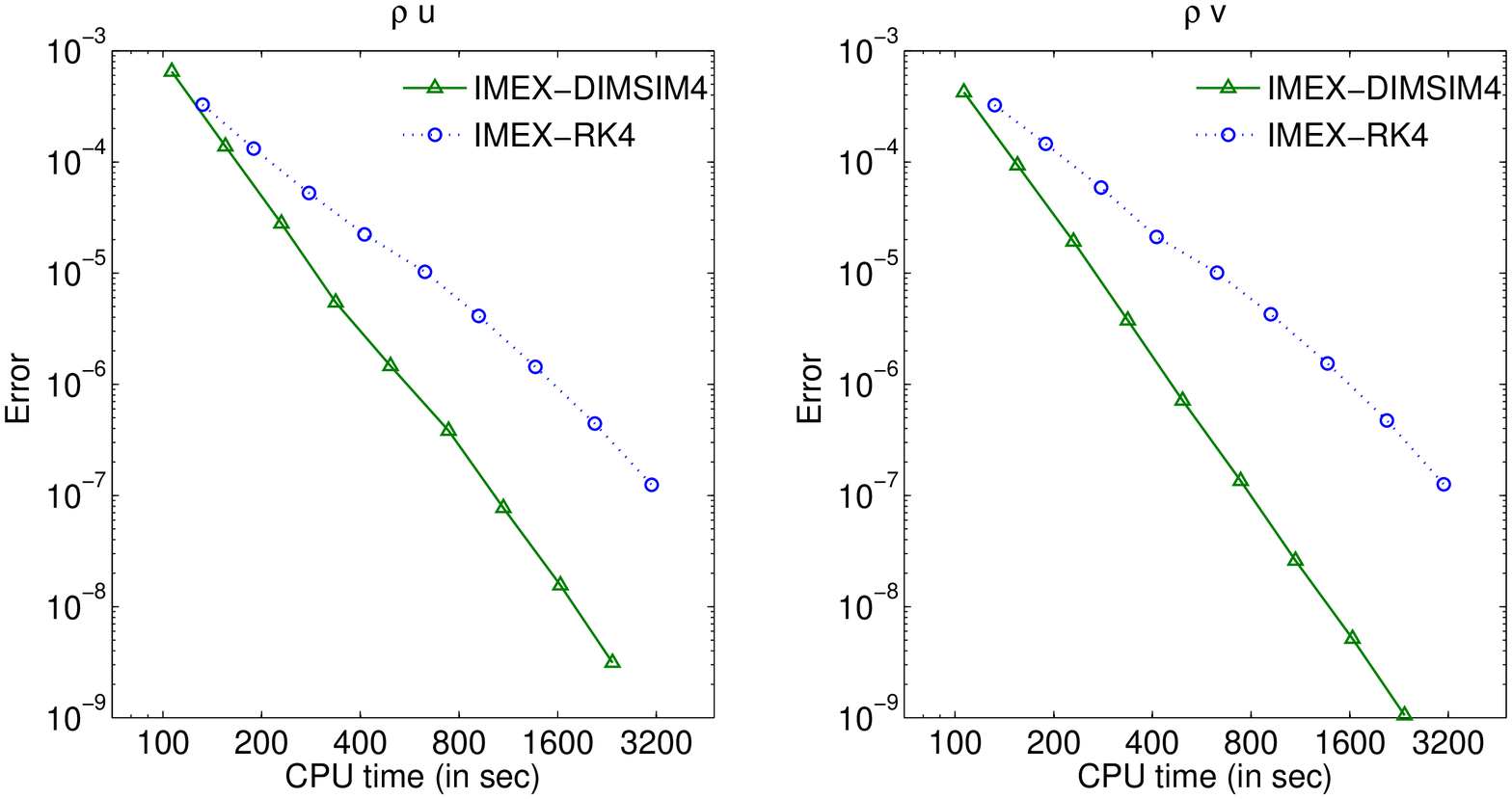}
  }
  }
\caption{Comparison of high order IMEX-DIMSIM and IMEX-RK results for the 2D rising bubble simulation \eqref{eqn:new-euler}. The integration time interval is $[0,200]$ sec. and is divided into $400,600,900,1350,2025,3037,4555,6832,10248$ equal time steps to obtain the points in the diagrams.  Temporal errors for all the variables \eqref{eqn:L2err} are computed for the solution at the final time.}
\label{fig:risingbubble2D}
\end{figure}

\begin{figure}
  \subfigure[Convergence diagram]{
    \centering{
      \includegraphics[width=\textwidth]{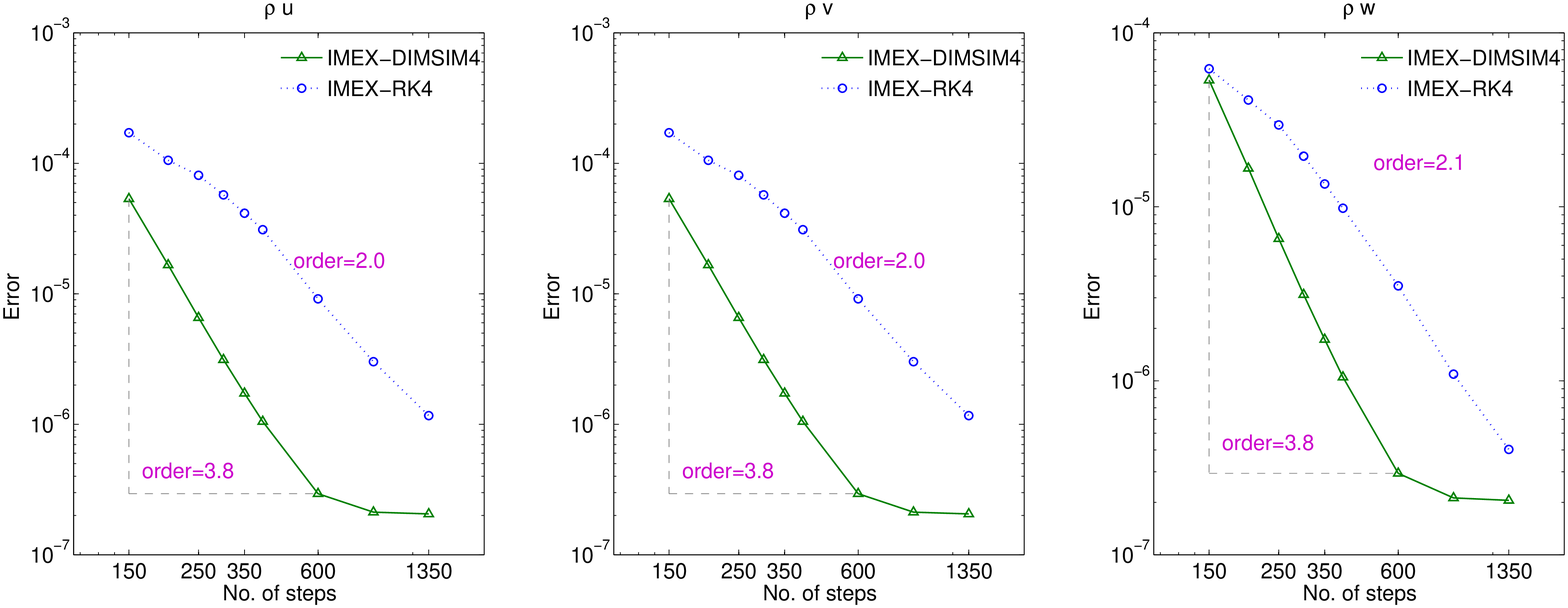}
    }
  }
  \subfigure[Work-precision diagram]{
  \centering{
  \includegraphics[width=\textwidth]{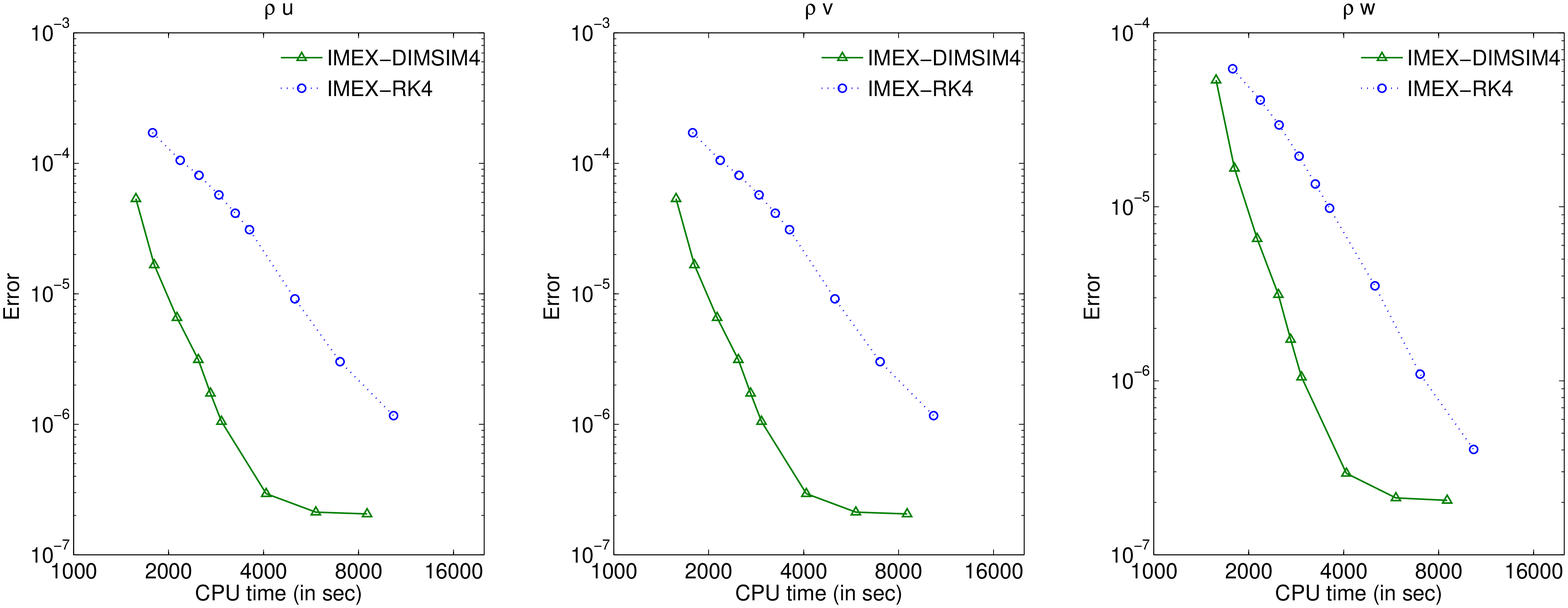}
  }
  }
\caption{Comparison of high order IMEX-DIMSIM and IMEX-RK results for the 3D rising bubble \eqref{eqn:new-euler}. The integration time interval is $[0,300]$ sec. and is divided into $150,200,250,300,350,400,600,900,1350$ equal time steps to obtain the points in the diagrams. Temporal errors for all the variables \eqref{eqn:L2err} are computed for the solution at the final time.
}
\label{fig:risingbubble3D}
\end{figure}

\begin{figure}
    \centering{
      \includegraphics[width=0.58\textwidth]{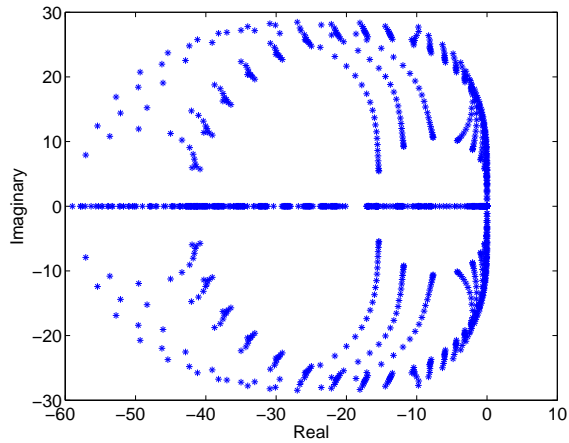}
    }
\caption{Plot of eigenvalues of the Jacobian for 2D rising bubble test problem. }
\label{fig:risingbubble2D_eigen}
\end{figure}

Figure \ref{fig:risingbubble2D} compares the convergence results and efficiency for the fourth-order IMEX-DIMSIM and IMEX-RK methods for  the 2D simulations.  As expected, the IMEX-DIMSIM reproduces the theoretical order of accuracy.
But the IMEX-RK scheme shows an obvious order reduction, which translates into a loss of computational efficiency. 
The 3D results are given in Figure \ref{fig:risingbubble3D}. 
The IMEX-RK method stills yields order reduction, less severely though. 
The error behavior of the IMEX-DIMSIM is somewhat irregular. 
It shows high order in the beginning, and then plateaus at the accuracy level $10^{-7}$ for a wide range of decreasing step sizes. 
The error plateau is likely due to the level of accuracy of the reference solution. 
However, even with this irregular behavior, the IMEX-DIMSIM is considerably more efficient than the IMEX-RK method. 
We have also tested large step sizes and found that the maximal allowable step size for IMEX-RK4 and IMEX-DIMSIM4 are both approximately equal to $1.0$ sec. 
This agrees with the prediction of the stability analysis in section \ref{sec:properties} which shows that the IMEX-DIMSIM has a good stability property. 
Furthermore, we notice that neither IMEX-RK5 nor IMEX-DIMSIM5 is suitable for this test problem because the maximal step sizes for them are restricted to values that are too small to make them competitive. 
Figure \ref{fig:risingbubble2D_eigen} shows that there are many eigenvalues of the Jacobian close to the imaginary axis, 
therefore a stability region covering a large part of imaginary axis is highly desirable. 
%
\section{Conclusions and future work}\label{sec:conclusions}

Multiscale problems in science and engineering are modeled by time-dependent systems of equations involving both stiff and nonstiff terms.
Implicit-explicit time stepping schemes perform an implicit integration only for the stiff components of the system, and thus combine the low cost of explicit methods with the favorable stability properties of implicit methods.

Many modern PDE solvers use high order spatial discretization schemes, e.g., the discontinuous Galerkin approach with high degree polynomials. Often the high order of spatial discretization is paired with a low order traditional time stepping scheme. It is therefore of considerable importance to develop high order time stepping algorithms that match the accuracy of the spatial discretization.

This paper addresses the need for high order implicit-explicit temporal discretizations in large scale applications.
We construct new fourth and fifth order IMEX DIMSIM schemes based on L-stable implicit components, and with the explicit components optimized such as to maximize the constrained stability regions.  The new methods have good stability properties and can take large step sizes for stiff problems.

Several test problems from different application areas  that can benefit from implicit-explicit integration are considered.
These problems are the two-dimensional Allen-Cahn and Burgers' equations with finite difference spatial discretizations, and two- and three-dimensional compressible Euler equations with discontinuous Galerkin space discretizations.
The performance of the new fourth and fifth order IMEX-DIMSIMs is compared against existing fourth and fifth order IMEX-RK methods.
In all cases the IMEX-DIMSIMs can use large step sizes  - similar to those taken by traditional implicit-explicit Runge-Kutta methods.  However, the high stage order enables our methods to avoid the order reduction that plagues classic IMEX-RK methods when applied to stiff systems or to problems with complex boundary conditions. In all cases IMEX-DIMSIMs are considerably more efficient than traditional IMEX-RK methods of the same order.


Typically multiscale flow simulations are carried out using fixed, predefined time steps. This is the approach
taken in this paper as well. On-going work by the first two authors focuses on the development of adaptive stepsize IMEX-GLM schemes.  

The high order IMEX-GLM schemes proposed herein are not only of interest to multiscale nonhydrostatic atmospheric simulations, 
but also to many other fields where large-scale multiscale simulations are carried out with high order spatial discretizations. 
IMEX-GLMs can prove especially useful in situations where IMEX-RK methods suffer from order reduction; specific examples include stiff systems of singular perturbation type or problems with challenging time-dependent boundary conditions.

\section*{Acknowledgements}
This work was supported by the National Science Foundation  through the awards NSF DMS--0915047, NSF CCF--0916493, NSF OCI--0904397, NSF CMMI--1130667, NSF CCF--1218454, AFOSR FA9550--12--1--0293--DEF, and AFOSR 12--2640--06.


\bibliographystyle{siam}
\bibliography{zhang}
\begin{sidewaystable}
\small
\begin{eqnarray*}
&& A = \left[~
\begin{array}{r r r r}
                   0 &                  0 &                 0 & 0  \\
   0.258897065974412 &                  0 &                 0 & 0  \\
   2.729801825357062 & -0.060004247312668 &                 0 & 0  \\
   0.951308318232761 &  0.614160494289040 & 0.422498793609078 & 0  \\
\end{array}
~\right] \\
&& B = \left[~
\begin{array}{r r r r}
   5.669708110906782 & -0.493235358869745 &  0.021475944586626 &  0.175951726795284 \\
   5.544708110906782 &  0.020653530019144 & -0.797968499857818 &  0.680943549709761 \\
   4.720814974705226 &  3.191226074825372 & -5.227438428178271 &  0.686166890688894 \\
   4.848863779632135 &  2.337640759837926 & -3.218585217497575 &  0.418013495315584 \\
\end{array}
~\right] \\
&& Q = \left[~
\begin{array}{r r r r r}
   1 &                  0 &                  0 &                  0 &                  0 \\
   1 &  0.074436267358921 &  0.055555555555556 &  0.006172839506173 &  0.000514403292181 \\
   1 & -2.003130911377728 &  0.242223637993112 &  0.052716285344531 &  0.008600849263247 \\
   1 & -0.987967606130879 &  0.013613972830935 &  0.038658018404147 &  0.017011414548385 \\
\end{array}  
~\right] \\   
&& \widehat{A} = \left[~
\begin{array}{r r r r}
   0.572816062482135 &                  0 &                  0 &                  0 \\
   0.294478591621391 &  0.572816062482135 &                  0 &                  0 \\
   3.754531024312379 & -0.446626145372372 &  0.572816062482135 &                  0 \\
  20.906355951077522 & -6.918033573971423 &  0.824272703722306 &  0.572816062482135 \\
\end{array}  
~\right] \\
&& \widehat{B} = \left[~
\begin{array}{r r r r}
   2.818382755109841 &  -0.107847984112942 &  1.213319973963157 & -0.548700992864529 \\
   3.266198817591976 & -1.885223345152593  & 3.830771904411522 & -1.797738883043436 \\
   3.774131970777119 & -3.469139895411032  & 5.100995462482731 & -4.672071998026633 \\
   1.800600620848989 &  6.203817506581311  &-13.407704583723200 & -5.034154872439978 \\
\end{array}  
~\right] \\
&& \widehat{Q} = \left[~
\begin{array}{r r r r r}
   1 & -0.572816062482135 &                  0 &                  0  &                 0 \\
   1 & -0.533961320770192 & -0.135383131938489 & -0.025650275076168 & -0.003021498328079 \\
   1 & -3.214054274755475 & -0.010779770975077 & -0.053097178648182 & -0.017299808772539 \\
   1 & -14.385411143310540 &  1.683679993026802 &  0.081422122041277 & -0.051803591005091 \\
\end{array}  
~\right] \\   
&& v = [~0.281364340879037 \quad -1.282889560784121 \quad  2.266595749735792 \quad -0.265070529830707~] \\
&& c = [~0 \quad 1/3 \quad 2/3 \quad 1~]
\end{eqnarray*}
\caption{\label{tab:order-4-implicit-coefficients}Coefficients of the IMEX-DIMSIM-4.}
\end{sidewaystable}
\begin{sidewaystable}
\small
\begin{eqnarray*}
&& A = \left[~
\begin{array}{r r r r r}
                   0 &                  0 &                  0 &                  0 & 0 \\
   0.380631951399918 &                  0 &                  0 &                  0 & 0 \\
  -0.723344119927179 &  0.934338548518619 &                  0 &                  0 & 0 \\
  -0.292421654731536 &  1.489386717103117 &  0.229042913082062 &                  0 & 0 \\
  10.333193352608074 &  0.200217292186561 &  0.841800685401247 & -0.148918889975160 & 0 \\
\end{array}
~\right] \\
&& B = \left[~
\begin{array}{r r r r r}
  -1.811278483713069 &  2.072219536433343 &  0.130011155311711 &  0.166279568600910 &  0.117403740739418 \\
  -1.724125705935292 &  1.629858425322231 &  1.038344488645044 & -0.796914875843534 &  0.396841233783945 \\
  -1.998394810009466 &  3.088356723470882 & -2.146707663207811 &  2.854109498231544 & -0.833722659704275 \\
  -1.361504766226497 &  0.334933035918415 &  2.154212895587752 &  0.353113262914561 & -1.482126886275562 \\
   5.091061924499312 &-29.458910962376240 & 55.143920860593482 &-43.440447985319850 &  3.112719239754878 \\
\end{array}
~\right] \\
&& Q = \left[~
\begin{array}{r r r r r r}
   1 &                  0 &                  0 &                  0 &                  0 &                  0 \\
   1 & -0.130631951399918 &  0.031250000000000 &  0.002604166666667 &  0.000162760416667 &  0.000008138020833 \\
   1 &  0.289005571408560 & -0.108584637129655 & -0.008364746307874 &  0.000170993363233 &  0.000108343335202 \\
   1 & -0.676007975453643 & -0.205618135816810 & -0.004861199044730 &  0.004533255151668 &  0.001138659940362 \\
   1 &-10.226292440220721 &  0.140734501734106 &  0.097068228416195 &  0.034078612640450 &  0.008071842745668 \\
\end{array}  
~\right] \\   
&& \widehat{A} = \left[~
\begin{array}{r r r r r}
   0.278053841136452 &                  0 &                  0 &                  0 &                  0 \\
   0.220452276182580 &  0.278053841136452 &                  0 &                  0 &                  0 \\
   2.294819895736366 & -0.602366708071285 &  0.278053841136452 &                  0 &                  0 \\
   5.054620901153854 & -1.529876218309763 &  0.097119141498823 &  0.278053841136452 &                  0 \\
   9.345167780108133 & -1.412133513099773 & -1.883401998517870 &  0.782533955446870 &  0.278053841136452 \\
\end{array}  
~\right] \\
&& \widehat{B} = \left[~
\begin{array}{r r r r r}
   6.044855283302179 & -2.020000467205476 &  0.032934533641225 &  0.593578985923315 & -0.226664851205853 \\
   5.853954219943505 & -1.072092372634326 & -1.839270544389963 &  2.410922952843391 & -0.899263047489796 \\
   6.004175007913425 & -2.014097375842605 &  0.610845429880394 & -0.963490004887004 & -0.405182760273902 \\
   6.002703177071046 & -2.556003283230891 &  3.151551366098853 & -5.493514217893924 &  0.448102618067392 \\
   4.481882795290198 &  2.672564354868939 & -1.413660973235832 & -8.058154793746990 &  0.909905877341711 \\
\end{array}  
~\right] \\
&& \widehat{Q} = \left[~
\begin{array}{r r r r r r}
   1 &  -0.278053841136452 &                  0 &                  0 &                  0 &                  0 \\
   1 &  -0.248506117319032 & -0.038263460284113 & -0.006085015868847 & -0.000561338127960 & -0.000037118138206 \\
   1 &  -1.470507028801533 &  0.136564756449595 &  0.004900562818504 & -0.001619958388074 & -0.000365640421568 \\
   1 &  -3.149917665479366 &  0.406619102975690 &  0.027778596315200 & -0.004406329750951 & -0.001692120959916 \\
   1 &  -6.110220065073812 &  0.929780069812273 &  0.087106493228110 & -0.016782586272280 & -0.008434321001423 \\
\end{array}  
~\right] \\   
&& v = [~ -0.079385465132435 \quad  0.554317572910577 \quad -1.569589549144155 \quad  2.332074592443682  \quad-0.237417151077669~] \\
&& c = [~0 \quad 1/4 \quad 1/2 \quad 3/4 \quad 1~]
\end{eqnarray*}
\caption{\label{tab:order-5-implicit-coefficients}Coefficients of the IMEX-DIMSIM-5.}
\end{sidewaystable}

\end{document}